\documentclass[12pt, times]{article}
\usepackage{graphicx}
\usepackage{amssymb,amsmath}
\def\C{\mathcal{C}}
\def\B{\mathcal{B}}

\def\G{\mathcal{G}}
\def\H{\mathcal{H}}
\def\K{\mathcal{K}}

\def\P{\mathcal{F}}

\def\T{\mathcal{T}}

\def\N{\mathcal{N}}
\def\M{\mathcal{M}}
\def\Aut{\mathrm{Aut}}
\def\aut{\mathrm{aut}}

\newcommand{\Tc}{\T_{\mathrm{C}}}

\newcommand{\nn}{\mbox{\small 1} \hspace{-0,30em} 1}
\newtheorem{defn}{Definition}
\newtheorem{thm}{Theorem}
\newtheorem{corol}{Corollary}

\newtheorem{propos}{Proposition}

\begin{document}
\title{ Counting unlabelled toroidal graphs with no $K_{3,3}$-subdivisions\footnote{With the partial support of NSERC (Canada)}}
\bigbreak
\author{Andrei Gagarin, Gilbert Labelle and Pierre Leroux \\[0.1in]
\small \it Laboratoire de Combinatoire et d'Informatique Math\'ematique (LaCIM),\\ 
\small \it Universit\'e du Qu\'ebec \`a Montr\'eal (UQAM), Montr\'eal, Qu\'ebec, CANADA, H3C 3P8\\[0.1in]
\small e-mail: \texttt{gagarin@lacim.uqam.ca}, \texttt{labelle.gilbert@uqam.ca} and \texttt{leroux.pierre@uqam.ca}}
\maketitle
\begin{abstract}
We provide a description of unlabelled enumeration techniques, with complete proofs, for graphs that can be canonically obtained by substituting 2-pole networks for the edges of core graphs. 
Using structure theorems for toroidal and projective-planar graphs containing no $K_{3,3}$-subdivisions, we apply these techniques to obtain their unlabelled enumeration.
\end{abstract}
\section{Introduction}
We are interested in non-planar (finite, simple) graphs that can be embedded on the torus or the projective plane.
By Kuratowski's theorem \cite{Kuratowski}, a graph $G$ is planar if and only if it contains no subdivision of $K_5$ or of $K_{3,3}$, and, by Wagner's theorem \cite{Wagner}, a graph $G$ is planar if and only if it has no minor isomorphic to $K_5$ or $K_{3,3}$.
%In this paper 
Here, we restrict our attention to the graphs with no $K_{3,3}$-subdivisions. Since $K_{3,3}$ is a 3-regular graph, it is possible to see that the graphs with no $K_{3,3}$-subdivisions are precisely the graphs with no $K_{3,3}$-minors. Therefore we may refer to them as to {\em $K_{3,3}$-free graphs}. Characterizations of $K_{3,3}$-free toroidal graphs in terms of forbidden minors and forbidden subdivisions are given by Gagarin, Myrvold and Chambers in \cite{CGM,GMClong}. 

In \cite{GLL} and \cite{GLL2}, we have established structure theorems characterizing the classes of non-planar $K_{3,3}$-free 2-connected projective-planar graphs (denoted by $\P$) and toroidal graphs (denoted by $\T$) in terms of a special substitutional operation $\G\uparrow\N$, where 2-pole networks of a given class $\N$ are substituted for the edges of graphs from a class $\G$. Moreover, we have used these structure theorems to enumerate the labelled graphs in $\P$ and in $\T$. 

In the present paper we concentrate on the more difficult problem of their unlabelled enumeration. Our approach is based on Walsh's method \cite{Timothy} for the unlabelled enumeration of three-connected and homeomorphically irreducible two-connected graphs, together with species theory techniques. The paper is organized as follows. The techniques and theorems are presented in Sections 2 to 5. The proofs of the main results appear in Section 6. In Section 2, we recall the definition of the composition operation $\G\uparrow\N$ and review the structure theorems of \cite{GLL} and \cite{GLL2} for 2-connected non-planar $K_{3,3}$-free projective-planar and toroidal graphs. 

In Section 3, we introduce what we call the \emph{Walsh index series} %$W_\G$ 
of classes of graphs and networks. These are similar to the cycle index series of graphs except that the edge cycles induced by graph automorphisms are also taken into account. We first give their basic properties, in particular, the enumerative formulas related to classes of structures of the form $\G\uparrow\N$, and then proceed with their formal definitions.

In section 4 we calculate the necessary Walsh index series for the projective-planar and toroidal cores.  These are the basic graphs into which planar networks are substituted to form the graphs in $\P$ and $\T$.  There is an infinite family $\H$ of toroidal cores, called \emph{toroidal crowns}, whose Walsh series involves matching polynomials for paths and cycles and an extension of the substitution formulas for $\G\uparrow\N$-structures to the case where $\G$ represents graph structures with two sorts of edges.

Some numerical results are given in Section 5.  In particular the number of unlabelled toroidal crowns and of toroidal cores according to the number of vertices and edges is given, up to 64 vertices.  However, the enumeration of unlabelled non-planar 2-connected $K_{3,3}$-free projective-planar or toroidal graphs requires the computation of some enumerative series for planar networks which are still unknown in general. Results for small sizes can be obtained by generating these networks by hand or by computer algorithms. Tables 2 and 3 cover unlabelled non-planar 2-connected $K_{3,3}$-free projective-planar graphs on up to 9 vertices and non-projective-planar toroidal graphs on up to 12 vertices. 
%We hope to extend these tables significantly in the near future.

Finally, in Section 6, we provide complete and detailed proofs of the basic enumerative formulas, in particular Theorems 3 and 4 below. This involves redefining the Walsh series $W_{\G}$ and the tilde generating series $\widetilde{\G}$, 
for the unlabelled enumeration of $\G$, in terms of labelled enumeration of associated structures consisting of pairs $(G,\sigma)$, where $G$ is in $\G$ and $\sigma$ is an automorphism of $G$, with weights, in a typical species theory way. In other words, we transform the desired formulas into equations between exponential generating functions of classes of labelled structures. These formulas are then much easier to establish.
%
%\section{Decomposition and structure theorems}
%
\section{Decomposition and structure theorems}
By convention, the graph $K_2$ is considered as a 2-connected (non-separable) graph in this paper. A \textit{$2$-pole network} (or simply a \textit{network}) is a connected graph $N$ with two distinguished vertices $0$ and $1$, such that the graph $N\cup 01$ is $2$-connected, where the notation $N\cup ab$ is used for the graph obtained from $N$ by adding the edge $ab$ if it is not already there. The vertices $0$ and $1$ are called the {\it poles} of $N$, and all the other vertices of $N$ are said to be \textit{internal}. 

We define an operator $\tau$ acting on $2$-pole networks, $N \mapsto \tau\cdot N$, which interchanges the poles $0$ and $1$. A class $\N$ of networks is called \textit{symmetric} if $N\in\N \Rightarrow \tau\cdot N\in\N$. 
\smallskip

\noindent
{\bf Definition.}
Let $\G$ be a class of graphs and $\N$ be a symmetric class of networks.  We denote by $\G\uparrow \N$ the class of pairs of graphs $(G,G_0)$, such that
\begin{enumerate}
\item
the graph $G_0$ is in $\G$ ($G_0$ is called the \emph{core}),
\item
the vertex set $V(G_0)$ is a subset of $V(G)$,
\item
there exists a family $\{N_e:e\in E(G_0)\}$ of networks in $\N$ (called the \emph{components}) such that the graph $G$ can be obtained from $G_0$ by substituting  $N_e$ for each edge $e\in E(G_0)$, identifying the poles of $N$ with the extremities of $e$ according to some orientation.
\end{enumerate}

An example of a $(\G \uparrow \N)$-structure $(G,G_0)$, with $\G=P_4$, the class of path-graphs of order 4, and $\N$ is the class of all networks, is given in Figure \ref{fig:exemple}.
\begin{figure}[h] \label{fig:exemple}
\begin{center} \includegraphics[height=2.2in]{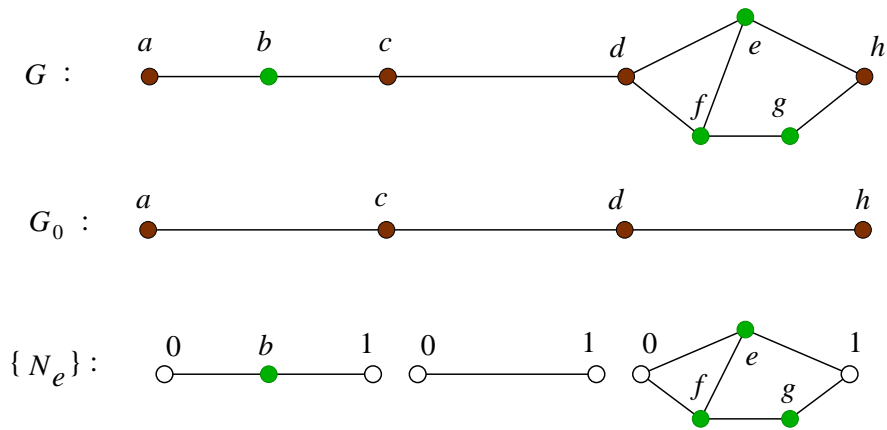}
\end{center}
\vspace{-4mm}
\caption{Example of a $(P_4 \uparrow \N)$-structure $(G,G_0)$}
\label{fig:exemple}
\end{figure}

We say that the composition $\G\uparrow \N$ is {\it canonical} if for any structure $(G,G_0)\in \G\uparrow \N$, the core  $G_0\in \G$ is uniquely determined by the graph $G$. In this case, we can identify $\G \uparrow \N$ with the class of resulting graphs $G$.

A network $N$ is called {\it strongly planar} if the graph $N\cup 01$ is planar. Denote by $\N_P$ the class of strongly planar networks.
In \cite{GLL} we prove the %uniqueness of the representation $\P=K_5\uparrow \N_P$ for 
following structure theorem for projective-planar graphs with no $K_{3,3}$-subdivisions whose proof is based on structural results of \cite{GK}.
\begin{thm}[\cite{GLL}]  The class $\P$ of $2$-connected, non-planar, $K_{3,3}$-free  and projective-planar graphs is characterized by the relation
\begin{equation} \label{eq:FK5flecheNP}
\P=K_5\uparrow \N_P,
\end{equation}
the composition being canonical.
\end{thm}

\begin{figure}[h]
\begin{center}
	\includegraphics[width=2.7in]{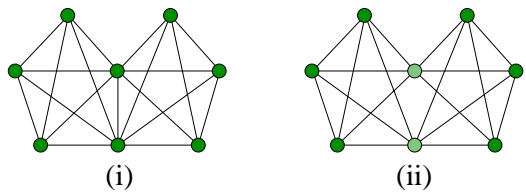}
\end{center}
\vspace{-5mm}
	\caption{(i) the graph $M$ \ \ \ (ii) the graph $M^*$}
	\label{fig:MetMetoile}
\end{figure}

\begin{figure}[h]
	\begin{center}
	\includegraphics*[width=3.7in]{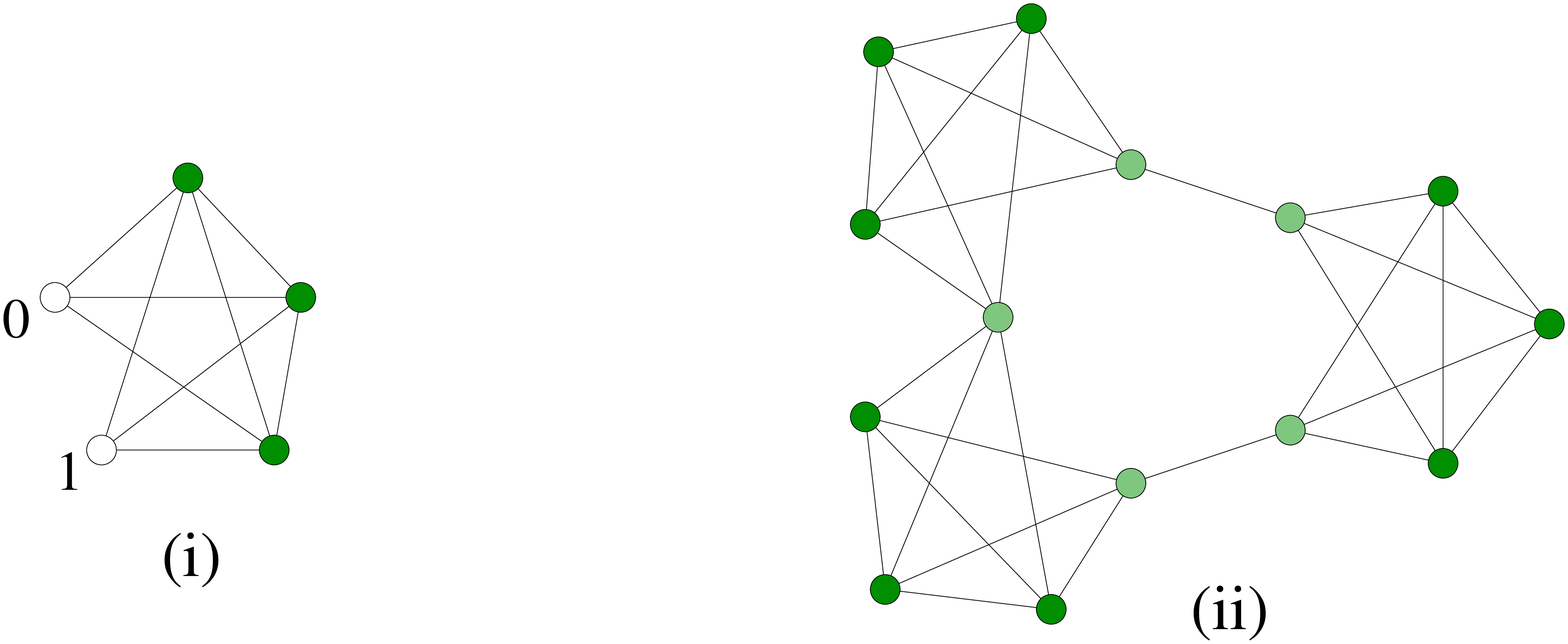}
	\end{center}
\vspace{-5mm}
	\caption{(i) a $K_5\backslash e$-network \ \ \ (ii) a toroidal crown}
	\label{fig:couronne}
\end{figure}

In order to describe a similar result of \cite{GLL2} for toroidal graphs, we need the following definitions. Given two disjoint $K_5$-graphs, the graph obtained by identifying an edge of one of the $K_5$'s with an edge of the other is called an {\it $M$-graph} (see Figure~\ref{fig:MetMetoile}(i)), and, when the edge of identification is deleted, an {\it $M^*$-graph} (see Figure~\ref{fig:MetMetoile}(ii)). A network obtained from $K_5$ by selecting two poles $0$ and $1$ among the vertices and by removing the edge $01$ is called a {\it $K_5\backslash e$-network} (see Figure~\ref{fig:couronne}(i)). 

Denote by $C_i$, $i\ge 3$, a cycle graph on $i$ vertices. A {\it toroidal crown} % of $K_5\backslash e$-networks 
$H$ is a graph obtained from a cycle $C_i$, $i\ge 3$, by substituting $K_5\backslash e$-networks for some edges of $C_i$ in such  a way that no two unsubstituted edges of $C_i$ are adjacent in $H$ (see Figure~\ref{fig:couronne}(ii)). Denote by $\H$ the class of toroidal crowns.

A {\it toroidal core} is defined as either $K_5$, an $M$-graph, an $M^*$-graph, or a toroidal crown. Denote by $\T_C$ the class of toroidal cores. In other words, $\T_C = K_5 + M + M^* + \H$. The following structure theorem is the basic result of \cite{GLL2}, obtained by refining the structural and algorithmic results of \cite{GK}.

\begin{thm}[\cite{GLL2}] \label{thm:thm2} 
The class $\T$ of $2$-connected non-planar $K_{3,3}$-free toroidal graphs is characterized by the relation
\begin{equation}\label{eq:TCflecheNP}
\T=\Tc\uparrow \N_P,
\end{equation}
the composition being canonical.
\end{thm}
%
%Theorems 1 and 2 imply that a projective-planar graph with no $K_{3,3}$-sub\-divisions is toroidal. However an arbitrary projective-planar graph can be non-toroidal. The characterizations of Theorems 1 and 2 can be used to detect toroidal or projective-planar graphs with no $K_{3,3}$'s in linear time. The implementation of this algorithm can be derived from \cite{GK} by using a breadth-first or depth-first search technique for the decomposition and by doing a linear-time planarity testing. The linear-time complexity of the detection follows from the linear-time complexity of the decomposition and from the fact that each vertex of the initial graph can appear in at most 7 different components.

%
\section{The Walsh index series of species of graphs and networks}
A \textit{species} is a class $\C$ of labelled combinatorial structures (for example, graphs) which is closed under isomorphism. Each combinatorial structure has an underlying set (for example, the vertex set of a graph), and any isomorphism is induced by a relabelling along a bijection between the underlying sets. Examples of species arise from classes of graphs which are closed under isomorphism and also from classes of networks, where the underlying set of a network consists of its internal vertices.

A species $\C$ is said to be \textit{weighted} if each structure $s$ of $\C$ is assigned a weight $w(s)$ taken from a commutative ring such that the weight function $s \longmapsto w(s)$ is invariant under isomorphism. For example, given a graph $G$, we can define the weight  $w_0(G)=y^m$, where $m=|E(G)|$ is the number of edges in $G$ and $y$ is a formal variable acting as an \textit{edge counter}. 

We write $\C=\C_w$ to denote a weighted species $\C$ with weight function $w$.
Given a weighted species $\C=\C_w$, we use two generating functions corresponding to the weighted enumeration of labelled and unlabelled $\C$-structures: the exponential generating function
\begin{equation} \label{eq:Cwx}
\C_w(x) = \sum_{n\ge 0} |\C[n]|_w\frac{x^n}{n!},
\end{equation}
where $\C[n]$ denotes the set of all $\C$-structures with underlying set $[n]=\{1,2,\ldots ,n\}$ (labelled structures) and where, for a weighted set $S$, the \emph{total weight} $|S|_w$ is defined by
\[ |S|_w=\sum_{s\in S}w(s),\]
and the ordinary generating function
\begin{equation} \label{eq:Cwtildex}
\widetilde{\C}_w(x)  = \sum_{n\ge 0} |\C[n]\slash\!\!\sim|_w x^n,
\end{equation}
where $\C[n]\slash \!\!\! \sim$ denotes the set of isomorphism classes of structures in $\C[n]$ (unlabelled structures), often called the \emph{tilde} generating function.
For example, for a species $\G$ of graphs weighted by the function $w_0(G)=y^m$, the generating functions (\ref{eq:Cwx}) and (\ref{eq:Cwtildex}) take the following form:
\begin{eqnarray}  
\G(x,y):=\G_{w_0}(x)=\sum_{n\ge 0}g_n(y)\frac{x^n}{n!},
\end{eqnarray}
with $g_n(y)=\sum_{m\ge 0}g_{n,m}y^m$, where $g_{n,m}$ is the number of graphs in $\G$ over the set $[n]$ of vertices and having $m$ edges, and
\begin{eqnarray}
\widetilde{\G}(x,y):=\widetilde{\G}_{w_0}(x)=\sum_{n\ge 0}\tilde{g}_n(y)x^n,
\end{eqnarray}
with $\tilde{g}_n(y)=\sum_{m\ge 0}\tilde{g}_{n,m}y^m$, where $\tilde{g}_{n,m}$ is the number of isomorphism classes of graphs in $\G$ having $n$ vertices and $m$ edges.

Similarly, for a species $\N$ of 2-pole networks weighted by the function $w_0(N)=y^m$, the generating functions (\ref{eq:Cwx}) and (\ref{eq:Cwtildex}) take the form:
\begin{eqnarray}
\N(x,y):=\N_{w_0}(x)=\sum_{n\ge 0}\nu_n(y)\frac{x^n}{n!},
\end{eqnarray}
with $\nu_n(y)=\sum_{m\ge 0}\nu_{n,m}y^m$, where $\nu_{n,m}$ is the number of networks in $\N$ over a set of $n$ internal vertices and having $m$ edges, and
\begin{eqnarray}
\widetilde{\N}(x,y):=\widetilde{\N}_{w_0}(x)=\sum_{n\ge 0}\tilde{\nu}_n(y)x^n,
\end{eqnarray}
with $\tilde{\nu}_n(y)=\sum_{m\ge 0}\tilde{\nu}_{n,m}y^m$, where $\tilde{\nu}_{n,m}$ is the number of unlabelled networks from $\N$ having $n$ internal vertices and $m$ edges. Note that unlabelled 2-pole networks are isomorphism classes of networks, where any isomorphism $\varphi : N \tilde{\longrightarrow} N^\prime$ is assumed to be pole preserving, i.e. $\varphi(0)=0$ and $\varphi(1)=1$. In particular, any automorphism of a network $N$ should be pole-preserving. 

We say that a species $\N$ of networks is \textit{symmetric} if the opposite network $\tau \cdot N$ of any network $N$ in $\N$ is also in $\N$, where the \emph{opposite} network $\tau \cdot N$ of a network $N$ is obtained by interchanging the poles 0 and 1 of $N$.
We also need to consider the subclass $\N_\tau$ of $\N$ consisting of  \textit{$\tau$-symmetric networks}, %in $\N$, 
that is, networks $N$ such that $\tau \cdot N$ is isomorphic to~$N$.

Suppose we are given a species $\G$ of graphs and a symmetric species $\N$ of networks. %such that the composition $\G \uparrow \N$ is canonical. %Then it is possible to compute the generating functions $(\G \uparrow \N)(x,y)$ and $(\G \uparrow \N)\tilde\ (x,y)$ for the species of $\G \uparrow \N$. 
In \cite{GLL} we have shown, following Walsh \cite{Tim}, how to compute the exponential generating function $(\G \uparrow \N)(x,y)$. %for the labelled graphs of $\G \uparrow \N$. 
In fact, we have 
\begin{equation} \label{eq:gflechenxy}
(\G \uparrow \N)(x,y)=\G(x,\N(x,y)).
\end{equation}
%

%It is convenient to extend the definition of $\G\uparrow\N$ to the case where the composition is not canonical. % when $\N$ is a symmetric species of networks. $\G\uparrow \N$ is then defined to be the species whose structures consist of graphs $G\in\G\uparrow\N$ together with a selected core $G_0\in\G$ of $G$ and induced components $N_1,N_2,\allowbreak\ldots, N_k$ (up to pole interchanges) in $\N$, and whose isomorphisms are the core-preserving graph isomorphisms. For example, if we take $\G=K_2$, then the $(K_2\uparrow \N_P)$-structures consist of graphs $G$ together with two selected (adjacent or not) vertices $a$ and $b$, such that the graph $G\cup ab$ is 2-connected and planar. Isomorphisms of $(K_2\uparrow \N_P)$-structures are graph isomorphisms preserving the two selected vertices. Recall that when the composition is canonical, the cores $G_0\in\G$ are uniquely determined so that they don't have to be specified. With this new definition of $\G\uparrow \N$, equation (\ref{eq:gflechenxy}) remains valid.
Note that isomorphisms of $(\G\uparrow\N)$-structures $(G,G_0)$ are defined as core-preserving graph isomorphisms.
In the present paper we describe and extend the method of Walsh \cite{Timothy} for computing the tilde generating function $(\G \uparrow \N)\tilde \ (x,y)$ for unlabelled $(\G\uparrow\N)$-structures. % Detailed proofs are given in Section 6. 
The method involves special cycle index series 
$$W_\G({\bf a}; {\bf b}; {\bf c}), %=W_\G(a_1, a_2,\ldots;\allowbreak b_1, b_2,\ldots;\allowbreak c_1, c_2,\ldots),$$$$
\ \ \ W^+_\N({\bf a}; {\bf b}; {\bf c})\ \ \  \mathrm{and}\ \ \  W^-_\N({\bf a}; {\bf b}; {\bf c})$$ 
in variables ${\bf a}=(a_1,a_2,\ldots)$, ${\bf b}=(b_1,b_2,\ldots)$ and ${\bf c}=(c_1,c_2,\ldots)$, which we call the \textit{Walsh index series}. These series are defined below (see Definitions\,\,1 and\,\,2). We first state several properties and theorems indicating their importance and application to the enumeration of graphs. %They have the following properties.
\begin{propos} Let $\G$ be a species of graphs. Then the following series identities hold:
\begin{eqnarray}
\G(x,y) & = & W_\G(x,0,0,\ldots; y,0,0,\ldots; 0,0,0,\ldots),
\label{form:7}\\
\widetilde{\G}(x,y) & = & W_\G(x,x^2,x^3,\ldots; y,y^2,y^3,\ldots; y,y^2,y^3,\ldots).
\label{form:8}
\end{eqnarray}
\end{propos}
\begin{propos} Let $\N$ be a species of networks. Then the following series identities hold:
\begin{eqnarray}
\widetilde{\N}(x,y) & = & W_{\N}^+(x,x^2,x^3,\ldots; y,y^2,y^3,\ldots; y,y^2,y^3,\ldots),\\
\widetilde{\N}_{\tau}(x,y) & = & W_{\N}^-(x,x^2,x^3,\ldots; y,y^2,y^3,\ldots; y,y^2,y^3,\ldots).
\end{eqnarray}
\end{propos}

Denote by 
\begin{eqnarray} \label{eq:W+Nk}
 (W^+_\N)_k=(W^+_\N)_k({\bf a}; {\bf b}; {\bf c}) = W^+_\N(a_k, a_{2k},\ldots; b_k, b_{2k},\ldots; c_k, c_{2k},\ldots)
\label{form:W+Nk}
\end{eqnarray}
and by 
\begin{eqnarray} \label{eq:W-Nk}
(W^-_\N)_k=(W^-_\N)_k({\bf a}; {\bf b}; {\bf c}) = W^-_\N(a_k, a_{2k},\ldots; b_k, b_{2k},\ldots; c_k, c_{2k},\ldots).
\label{form:W-Nk}
\end{eqnarray}
Then we have the following basic theorem.
\begin{thm} Let $\G$ be a species of graphs and $\N$ be a symmetric species of networks.%, such that the composition $\G \uparrow \N$ is canonical. 
Then the Walsh index series of the species $\G\uparrow\N$ is given by
\begin{eqnarray}
W_{\G \uparrow \N}({\bf a}; {\bf b}; {\bf c}) & = & W_\G(a_1,a_2,\ldots; (W^+_\N)_1,(W^+_\N)_2,\ldots; (W^-_\N)_1,(W^-_\N)_2,\ldots).
\end{eqnarray}

\end{thm}

As a corollary, we obtain the generating function $(\G\uparrow\N)\tilde\ (x,y)$. This is stated in the following theorem. 
\begin{thm} Given a species $\G$ of graphs and a symmetric species $\N$ of networks, %such that the composition $\G \uparrow \N$ is canonical, 
the generating function $(\G\uparrow\N)\tilde\ (x,y)$ of unlabelled $(\G\uparrow\N)$-structures is given by
\begin{eqnarray}
(\G \uparrow \N)\tilde \ (x,y) = W_\G(x,x^2,\ldots; \widetilde{\N}(x,y),\widetilde{\N}(x^2,y^2),\ldots; \widetilde{\N}_\tau(x,y),\widetilde{\N}_\tau(x^2,y^2),\ldots).
\end{eqnarray}
\end{thm}

\begin{figure}[h]
\begin{center}
	\includegraphics[width=4.2in]{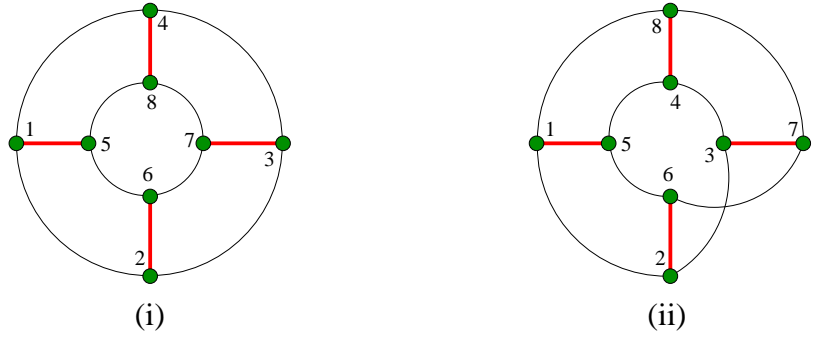}
\end{center}
\vspace{-5mm}
	\caption{(i) Cylindrical edge cycle, (ii) M\"obius edge cycle.}
	\label{fig:cyletmob}
\end{figure}

We now proceed to the formal definition of the Walsh index series $W_\G({\bf a}; {\bf b}; {\bf c})$ for a species $\G$ of graphs.
Let $G=(V(G),E(G))$ be a graph in $\G$. A permutation $\sigma$ of $V(G)$ which is an automorphism of the graph $G$, induces a permutation $\sigma^{(2)}$ of the set $E(G)$ of edges whose cycles are of two possible sorts: if $c$ is a cycle of $\sigma^{(2)}$ of length $l$, then either $\sigma^l(a)=a$ and $\sigma^l(b)=b$ for each edge $e=ab$ of $c$, in which case $c$ is called an orientation preserving (or \textit{cylindrical}) edge cycle, or else $\sigma^l(a)=b$ and $\sigma^l(b)=a$ for each edge $e=ab$ of $c$, in which case $c$ is called an orientation reversing (or \textit{M\"obius}) edge cycle. For example, the automorphism $\sigma=(1,2,3,4)(5,6,7,8)$ of the graph of Figure~\ref{fig:cyletmob}(i) induces the cylindrical edge cycle 
%$(\{1,5\},\{2,6\},\{3,7\},\{4,8\})$, 
$(15,26,37,48)$ and the automorphism $\sigma=(1,2,3,4,5,6,7,8)$ of the graph of Figure~\ref{fig:cyletmob}(ii) induces the M\"obius edge cycle %$(\{1,5\},\{2,6\},\{3,7\},\allowbreak\{4,8\})$. 
$(15,26,37,48)$.

Let $\Aut(G)$ denote the set of automorphisms of $G$. For $\sigma\in\Aut(G)$, denote by $\sigma_k$ the number of cycles of length $k$ of the permutation $\sigma$, by $\mathrm{cyl}_k(G,\sigma)$ the number of cylindrical edge cycles of length $k$ of $G$ induced by $\sigma$, and by $\mathrm{m\ddot{o}b}_k(G,\sigma)$ the number of M\"obius edge cycles of length $k$ of $G$ induced by $\sigma$.
Given a graph $G\in \G$ and an automorphism $\sigma$ of $G$, we define the weight $w(G,\sigma)$ of such a structure as the following cycle index monomial:
\begin{eqnarray}
w(G,\sigma)=a_1^{\sigma_1}a_2^{\sigma_2}\cdots
b_1^{\mathrm{cyl}_1(G,\sigma)}b_2^{\mathrm{cyl}_2(G,\sigma)}\cdots c_1^{\mathrm{m\ddot{o}b}_1(G,\sigma)}c_2^{\mathrm{m\ddot{o}b}_2(G,\sigma)}\cdots\,.
\label{form:weightGsigma}
\end{eqnarray}

\begin{defn}
The Walsh index series $W_\G({\bf a}; {\bf b}; {\bf c})$ of a species $\G$ of graphs is defined as 
\begin{eqnarray}
W_\G({\bf a}; {\bf b}; {\bf c}) = \sum_{G\underline{\in} \mathrm{Typ}(\G)} \frac{1}{|\mathrm{Aut}(G)|} \sum_{\sigma\in \mathrm{Aut}(G)} w(G,\sigma),
\label{form:Wabc}
\end{eqnarray}
where the notation $G\underline{\in} \mathrm{Typ}(\G)$ means that the summation should be taken over a set of representatives $G$ of the isomorphism classes of graphs in $\G$.%, and $\mathrm{Aut}(G)$ denotes the automorphism group of $G$.
\end{defn}

We now define the Walsh index series $W^+_\N$ and $W^-_\N$ of a species $\N$ of 2-pole networks.
Let $\sigma\in S[U]$ be a permutation of the underlying set $U$ of a 2-pole network $N$. We can extend $\sigma$ to two permutations on $U\cup \{ 0,1 \}$, $\sigma^+=(0)(1)\:\!\sigma$
%\sigma\cup\{\mathrm{id}_{\{0,1\}} \}$ and $\sigma^{-}=\sigma\cup\{ \tau \}$, where $\tau$ is the transposition $(0,1)$. 
and $\sigma^{-}=(0,1)\:\!\sigma$. If $N$ is a $2$-pole network of $\N$ with internal vertices $U$, denote by $\hat{N}$ the corresponding graph on $U\cup \{ 0,1 \}$. Then we denote by 
\begin{eqnarray}
\mathrm{Aut}^+(N)=\{ \sigma\in S[U]\ |\ \sigma^+\in \mathrm{Aut}(\hat{N}) \}
\label{form:A+}
\end{eqnarray}
and 
\begin{eqnarray}
\mathrm{Aut}^-(N)=\{ \sigma\in S[U]\ |\ \sigma^-\in \mathrm{Aut}(\hat{N}) \}.
\label{form:A-}
\end{eqnarray}
Remark that $\Aut^+(N)=\Aut(N)$.

For $N\in \N[U]$ and $\sigma\in \Aut^+(N)$, we assign the weight
\begin{eqnarray}
w(N,\sigma)=\frac{w(\hat{N},\sigma^+)}{a_1^2},
\label{form:weightN+} 
\end{eqnarray}
and for $N\in \N[U]$ and $\sigma\in \Aut^-(N)$, we set 
\begin{eqnarray}
w(N,\sigma)=\frac{w(\hat{N},\sigma^-)}{a_2}.
\label{form:weightN-}
\end{eqnarray}
In other words, only the internal vertex cycles are accounted for. %Recall that these weights are monomial functions of ({\bf a}; {\bf b}; {\bf c}).
\begin{defn}
For a species of networks $\N$, we define the Walsh index series $W^+_\N$ and $W^-_\N$ by the formulas
\begin{eqnarray} 
W^+_\N({\bf a}; {\bf b}; {\bf c})= \sum_{N\underline{\in} \mathrm{Typ}(\N)} \frac{1}{|\mathrm{Aut}^+(N)|} \sum_{\sigma\in \mathrm{Aut}^+(N)} w(N,\sigma),
\end{eqnarray}
and
\begin{eqnarray} 
W^-_\N({\bf a}; {\bf b}; {\bf c})= \sum_{N\underline{\in} \mathrm{Typ}(\N_\tau)} \frac{1}{|\mathrm{Aut}^-(N)|} \sum_{\sigma\in \mathrm{Aut}^-(N)} w(N,\sigma).
\end{eqnarray}
\end{defn}

Let $\B$ be a species of $2$-connected graphs containing $K_2$. Denote by $\B_{0,1}$ the species of 2-pole networks that are obtained from the graphs of $\B$ by removing an edge and by relabelling the extremities  $0$ and $1$ in one or two possible ways. For example, for the class $K_5$ of complete graph with five vertices, we obtain the class $(K_5)_{0,1}$ of $K_5\backslash e$-networks. Also denote by $\N_{\B}$ the class of networks obtained by taking all networks in $\B_{0,1}$ except the trivial network $\nn$ consisting of two isolated vertices $0$ and $1$, and networks of the form $N\cup 01$, where $N\in \B_{0,1}$. Note that for the labelled enumeration, we have (see \cite{GLL} and the proof of Proposition 3 in Section 6)
\begin{equation} \label{eq:B01xy}
x^2\B_{0,1} (x,y)= 2\frac{\partial}{\partial y}\B(x,y)
\end{equation}
and
\begin{equation} \label{eq:NBxy}
\N_B(x,y) = (1+y)\B_{0,1} (x,y) - 1.
\end{equation}

\begin{propos}[\cite{Timothy}] Let $\B$ be a species of $2$-connected graphs containing $K_2$. Then the Walsh index series of the associated species of networks $\B_{0,1}$ and $\N_\B$ can be computed as follows:
\begin{eqnarray}
W^+_{\B_{0,1}}({\bf a}; {\bf b}; {\bf c})=\frac{2}{a_1^2}\frac{\partial}{\partial b_1}W_{\B}({\bf a}; {\bf b}; {\bf c}),
\label{eq:WB01+}
\end{eqnarray}
\begin{eqnarray}
W^-_{\B_{0,1}}({\bf a}; {\bf b}; {\bf c})=\frac{2}{a_2}\frac{\partial}{\partial c_1}W_{\B}({\bf a}; {\bf b}; {\bf c}),
\label{form:19}
\end{eqnarray}
\begin{eqnarray}
W^+_{\N_\B}({\bf a}; {\bf b}; {\bf c})=(1+b_1)W^+_{\B_{0,1}}({\bf a}; {\bf b}; {\bf c}) - 1
\label{form:20}
\end{eqnarray}
and
\begin{eqnarray}
W^-_{\N_\B}({\bf a}; {\bf b}; {\bf c})=(1+c_1)W^-_{\B_{0,1}}({\bf a}; {\bf b}; {\bf c}) - 1.
\label{form:21}
\end{eqnarray}
\end{propos}

\section{Walsh index series for toroidal cores}
Here we calculate the necessary Walsh index series for the toroidal and projective-planar cores of Theorems\,1 and\,2. These index series will be used later to count the unlabelled toroidal and projective-planar graphs with no $K_{3,3}$'s by using Theorem 4. 

First, let us consider the projective-planar and toroidal core graph $K_5$. Its Walsh index series is given in the following proposition.

\begin{propos}[\cite{Timothy}] 
\begin{eqnarray}
W_{K_5}({\bf a}; {\bf b}; {\bf c})=W(a_1,a_2,\ldots;b_1,b_2,\ldots;c_1,c_2,\ldots)\nonumber \\
=\frac{1}{|\Aut(K_5)|}\sum_{\sigma\in \Aut(K_5)} w(K_5,\sigma)=\frac{1}{5!}\sum_{\sigma \in S_5} w(K_5,\sigma)\nonumber \\
=\frac{1}{5!} ( a_1^5b_1^{10} + 10a_1^3a_2b_1^3b_2^3c_1 + 15a_1a_2^2b_2^4c_1^2 + 20a_1^2a_3b_1b_3^3\nonumber \\
+ 20a_2a_3b_3b_6c_1 + 30a_1a_4b_4^2c_2 + 24a_5b_5^2).
\end{eqnarray}
\end{propos} 

The corresponding $K_5\backslash e$-network has the following Walsh index series, using (\ref{eq:WB01+}):

\begin{corol}
\begin{eqnarray}
W^+_{K_5\backslash e}({\bf a}; {\bf b}; {\bf c})=
\frac{2}{a_1^2}\frac{\partial W_{K_5}({\bf a}; {\bf b}; {\bf c})}{\partial b_1}=
\frac{1}{3!}(a_1^3b_1^9 + 3a_1a_2b_1^2b_2^3c_1 + 2a_3b_3^3)
\label{form:w+k5}
\end{eqnarray}
and
\begin{eqnarray}
W^-_{K_5\backslash e}({\bf a}; {\bf b}; {\bf c})=
\frac{2}{a_2}\frac{\partial W_{K_5}({\bf a}; {\bf b}; {\bf c})}{\partial c_1}=
\frac{1}{3!}(a_1^3b_1^3b_2^3 + 3a_1a_2b_2^4c_1 + 2a_3b_3b_6).
\label{form:w-k5}
\end{eqnarray}
\end{corol}

\begin{propos} 
The toroidal core graph $M$ has the following Walsh index series
\begin{eqnarray}
W_M({\bf a}; {\bf b}; {\bf c})=\frac{1}{|\Aut(M)|}\sum_{\sigma\in \Aut(M)} w(M,\sigma)=
\frac{1}{144} [ a_1^8b_1^{19} + a_1^6a_2b_1^6b_2^6c_1\nonumber \\
+ 6(a_1^6a_2b_1^{12}b_2^3c_1 + a_1^4a_2^2b_1^3b_2^7c_1^2) + 9(a_1^4a_2^2b_1^5b_2^6c_1^2 + a_1^2a_2^3b_2^8c_1^3)\nonumber \\
+ 6(a_1^2a_2^3b_1b_2^9 + a_2^4b_2^9c_1) + 4(a_1^5a_3b_1^{10}b_3^3 + a_1^3a_2a_3b_1^3b_2^3b_3b_6c_1)\nonumber \\
+ 12(a_1^3a_2a_3b_1^3b_2^3b_3^3c_1 + a_1a_2^2a_3b_2^4b_3b_6c_1^2) + 4(a_1^2a_3^2b_1b_3^6 + a_2a_3^2b_3^2b_6^2c_1)\nonumber \\
+ 18(a_1^2a_2a_4b_1b_2^2b_4^3c_2 + a_2^2a_4b_2^2b_4^3c_1c_2) + 12(a_1^2a_6b_1b_6^3 + a_2a_6b_6^3c_1) ],
\end{eqnarray} 

and the toroidal core graph $M^*$ has the Walsh index series
\begin{eqnarray}
W_{M^*}({\bf a}; {\bf b}; {\bf c})=\frac{1}{|\Aut(M^*)|}\sum_{\sigma\in \Aut(M^*)} w(M^*,\sigma)=
\frac{1}{144} [ a_1^8b_1^{18} + a_1^6a_2b_1^6b_2^6\nonumber \\
+ 6(a_1^6a_2b_1^{11}b_2^3c_1 + a_1^4a_2^2b_1^3b_2^7c_1) + 9(a_1^4a_2^2b_1^4b_2^6c_1^2 + a_1^2a_2^3b_2^8c_1^2)\nonumber \\
+ 6(a_1^2a_2^3b_2^9 + a_2^4b_2^9) + 4(a_1^5a_3b_1^9b_3^3 + a_1^3a_2a_3b_1^3b_2^3b_3b_6)\nonumber \\
+ 12(a_1^3a_2a_3b_1^2b_2^3b_3^3c_1 + a_1a_2^2a_3b_2^4b_3b_6c_1) + 4(a_1^2a_3^2b_3^6 + a_2a_3^2b_3^2b_6^2)\nonumber \\
+ 18(a_1^2a_2a_4b_2^2b_4^3c_2 + a_2^2a_4b_2^2b_4^3c_2) + 12(a_1^2a_6b_6^3 + a_2a_6b_6^3) ].
\end{eqnarray}
\end{propos}

To obtain the Walsh index series $W_\H$ of the class $\H$ of toroidal crowns, we will use a variant of Theorem 3 where the $K_5\backslash e$-networks are substituted 
into some edges of a cycle $C_n\ (n\ge 3$) selected in such a way that no pair of unsubstituted edges are adjacent. In other words, the unsubstituted edges should form a matching of $C_n$.

Recall that a \textit{matching} $\mu$ of a finite graph $G$ is a set of pair-wise disjoint edges of $G$. We introduce the \textit{homogeneous matching polynomial} of $G$ as
$$M_G(y,z) = \sum_{\mu\in\M(G)} y^{|\mu|}z^{m-|\mu|},$$
where $\M(G)$ is the set of matchings of $G$ and $m=|E(G)|$. In particular, we have the homogeneous matching polynomials $U_n(y,z) = M_{P_n}(y,z)$ and $T_n(y,z) = M_{C_n}(y,z)$, where $P_n$ denotes the path graph and $C_n$, the cycle graph over the set of vertices $V=\{1,2,\ldots ,n\}$ (see \cite{GLL2}, \cite{Godsil}). They satisfy the recurrence relations
\begin{equation} \label{eq:Undeyzrec}
U_n(y,z) = yz U_{n-2}(y,z)+zU_{n-1}(y,z),
\end{equation}
\begin{equation} \label{eq:Tndeyzrec}
T_n(y,z) = y z^2 U_{n-2}(y,z)+zU_{n}(y,z),
\end{equation}
for $n\ge 3$, with the initial values $U_1(y,z)=1$, $U_2(y,z)=y+z$ and $T_1(y,z)=z$, $T_2(y,z)=2yz+z^2$. A useful convention is to also set $U_0(y,z)=1/z$. Notice that (\ref{eq:Undeyzrec}) and (\ref{eq:Tndeyzrec}) are then also valid for $n=2$.
These polynomials can also be computed by using their generating functions:
\begin{propos}
The ordinary generating function for the homogeneous matching polynomials $U_n(y,z)$ and $T_n(y,z)$ of paths and cycles, respectively, are given by
%\begin{eqnarray}
%\sum_{n\ge 1} U_n(y,z) x^n = \frac{x+x^2y}{1-xz-x^2yz}
%\end{eqnarray}
\begin{eqnarray}
\sum_{n\ge 0} U_n(y,z) x^n = \frac{1}{(1-xz-x^2yz)z}
\end{eqnarray}
and
\begin{eqnarray}
\sum_{n\ge 1} T_n(y,z) x^n = \frac{xz+2x^2yz}{1-xz-x^2yz}.
\end{eqnarray}
\end{propos}

Given a graph $G=(V,E)$ with a matching $\mu$, the edges are partitioned into two sorts, say $Y$ and $Z$, depending on their membership in $\mu$ (sort $Y$) or not (sort $Z$). Let $\G^\mathrm{m}$ denote the species of matched graphs, that is of pairs $(G,\mu)$, where $G$ is a graph in $\G$ and $\mu$ is a matching of $G$. Isomorphisms of matched graphs are matching-preserving graph isomorphisms. Let $\N$ be a species of 2-pole networks. We denote $\G^\mathrm{m}\uparrow_Z \N$ the class of graphs obtained by substituting networks of $\N$ into the edges of sort $Z$, i.e. the edges not belonging to the matching $\mu$ of the graph $G$. Our goal is to compute the Walsh index series $W_{\G^\mathrm{m}\uparrow_Z \N}$ of a class of graphs of the form $\G^\mathrm{m}\uparrow_Z \N$. 

In order to do this, we introduce an extension $W_\G^\mathrm{m}$ of the Walsh index series to the matched graphs as follows. Notice that an automorphism $\sigma$ of a matched graph $(G,\mu)$ is an automorphism $\sigma$ of $G$ which leaves the matching fixed, i.e. $\sigma(\mu) = \mu$. It will induce cylindrical and M\"obius edge cycles of sort $Y$, counted by the variables  $b_1,b_2,\ldots$ and $c_1,c_2,\ldots$, respectively, and cylindrical and M\"obius edge cycles of sort $Z$, counted by the variables $\beta_1,\beta_2,\ldots$ and $\gamma_1,\gamma_2,\ldots$, respectively.

Let $w^\mu(\sigma)$ denote the following cycle index monomial:
$$w^\mu(\sigma) = \prod_{k}a_k^{\sigma_k} \prod_{k}b_k^{\mathrm{cyl}_{k,Y}(\sigma)} \prod_{k}c_k^{\mathrm{m\ddot{o}b}_{k,Y}(\sigma)} \prod_{k}\beta_k^{\mathrm{cyl}_{k,Z}(\sigma)} \prod_{k}\gamma_k^{\mathrm{m\ddot{o}b}_{k,Z}(\sigma)},$$
where $\mathrm{cyl}_{k,Y}(\sigma), \mathrm{m\ddot{o}b}_{k,Y}(\sigma), \mathrm{cyl}_{k,Z}(\sigma)$ and  $\mathrm{m\ddot{o}b}_{k,Z}(\sigma)$ denote the number of cylindrical and M\"obius edge cycles of length $k$ and of sort $Y$ or $Z$, respectively. Then we set
$$W_\G^\mathrm{m}(\mathbf{a},\mathbf{b},\mathbf{c},\beta,\gamma)=\sum_{G\underline{\in}\mathrm{Typ}(\G)} \frac{1}{|\Aut(G)|}\sum_{\sigma\in\Aut(G)} \sum_{\mu\in \mathrm{Fix}^\mathrm{m}(\sigma)} w^\mu(\sigma),$$
where $\mathrm{Fix}^\mathrm{m}(\sigma)$ denotes the set of matchings $\mu$ of $G$ which are fixed by $\sigma$.

The following result can be seen as a corollary of Theorem 3. Recall the plethystic notation (\ref{form:W+Nk}) and (\ref{form:W-Nk}) for $(W_\N^+)_k$ and $(W_\N^-)_k$, respectively.
\begin{propos}
The Walsh index series $W_{\G^\mathrm{m} \uparrow_Z \N} (\mathbf{a};\mathbf{b};\mathbf{c})$ can be obtained from the extended Walsh series $W_\G^\mathrm{m}(\mathbf{a};\mathbf{b};\mathbf{c};\beta;\gamma)$ by performing the substitutions $\beta_k := (W_\N^+)_k(\mathbf{a};\mathbf{b};\mathbf{c})$, $\gamma_k := (W_\N^-)_k(\mathbf{a};\mathbf{b};\mathbf{c})$. In other words, 
$$W_{\G^\mathrm{m} \uparrow_Z \N} (\mathbf{a};\mathbf{b};\mathbf{c}) = W_\G^\mathrm{m}(\mathbf{a};\mathbf{b};\mathbf{c};(W_\N^+)_1,(W_\N^+)_2,\ldots;(W_\N^-)_1,(W_\N^-)_2,\ldots).$$
\end{propos}

Our main application of Proposition 7 is in the case where $\G=C=\sum_{n\ge 3}C_n$ and $\N=K_5\backslash e$. Indeed, by definition, for the species $\H$ of toroidal crowns, we have
$$\H=C^\mathrm{m} \uparrow_Z K_5\backslash e,$$
so that
\begin{eqnarray}
W_\H (\mathbf{a};\mathbf{b};\mathbf{c}) = W_C^\mathrm{m}(\mathbf{a};\mathbf{b};\mathbf{c};(W_{K_5\backslash e}^+)_1,(W_{K_5\backslash e}^+)_2,\ldots;(W_{K_5\backslash e}^-)_1,(W_{K_5\backslash e}^-)_2,\ldots).
\label{form:WH}
\end{eqnarray}
Thus we need to compute the extended Walsh series $W_{C_n}^\mathrm{m}(\mathbf{a};\mathbf{b};\mathbf{c};\beta;\gamma)$ for matched cycles of size $n$, $n\ge 3$. But first, as an example, we give the extended Walsh series $W_{P_n}^\mathrm{m}$ for matched paths of size $n\ge 1$.

\begin{propos} The extended Walsh series $W_{P_n}^\mathrm{m}(\mathbf{a};\mathbf{b};\mathbf{c};\beta;\gamma)$ of matched paths of size $n\ge 1$ is given by
\begin{eqnarray}
W_{P_n}^\mathrm{m}({\bf a}; {\bf b}; {\bf c}; \beta; \gamma)
=\frac{1}{2}a_1^nU_n(b_1,\beta_1) + 
\frac{1}{2}\left\{ \begin{array}{ll} 
a_1a_2^{\frac{n-1}{2}}\beta_2U_{\frac{n-1}{2}}(b_2,\beta_2), &\mbox{\emph \,$n$ }\mathrm{odd} \\ 
a_2^{\frac{n}{2}}(\gamma_1U_{\frac{n}{2}}(b_2,\beta_2) + c_1\beta_2U_{\frac{n-2}{2}}(b_2,\beta_2)),  &\mbox{\emph \,$n$ }\mathrm{even} 
\end{array}
\right .
\label{form:bcp}
\end{eqnarray}
\end{propos} 

\noindent \textbf{Proof.} Assume that the vertex set of $P_n$ is $[n]=\{1,2,\ldots,n\}$ and that the edges are of the form $\{i,i+1\}$ for $i=1,2,\ldots,n-1$. Then $\Aut(P_n)=\{\mathrm{id},\tau\}$, where $\mathrm{id}$ is the identity mapping of $[n]$ and $\tau$ is the reflection $(1,n)(2,n-1)\ldots$. For the identity mapping, any matching $\mu$ is left fixed and the edges of sort $Y$ (i.e. in $\mu$) become cylindrical cycles of length $1$, counted by $b_1$, and the edges of sort $Z$ give rise to cycles counted by $\beta_1$. This gives the first term $\frac{1}{2}a_1^nU_n(b_1,\beta_1)$. The second term corresponds to $\tau$-symmetric matchings. These are entirely determined by their restriction on the first half $\{1,2,\ldots,\lfloor\frac{n}{2}\rfloor\}$, and special attention has to be given to the parity of $n$. Details are left to the reader.
\hfill\rule{2mm}{2mm}

\medskip
%Cycle index of the oriented and non-oriented cycles:
%\begin{eqnarray}
%Z(C_n^\mathrm{or})=\frac{1}{n}\sum_{d|n}\varphi(d)x_d^{\frac{n}{d}}=\frac{1}{n}\sum_{d|n}\varphi(\frac{n}{d})x_\frac{n}{d}^{d},
%\end{eqnarray}

%\begin{eqnarray}
%Z(C_n^\mathrm{nor})=\frac{1}{2}Z_{C_n} + 
%\frac{1}{2}\left\{ \begin{array}{ll} 
%x_1x_2^{\frac{n-1}{2}}, & \mbox{\emph \,$n$ odd} \\ 
%\frac{1}{2}(x_2^{\frac{n}{2}} + x_1^2x_2^\frac{n-2}{2}), & \mbox{\emph \,$n$ even} 
%\end{array}.
%\right.
%\end{eqnarray}
We now move on to the matched cycles $C_n^\mathrm{m}$, $n\ge 3$. First, we give the Walsh index series for cycles $C_n$, $n\ge 3$, which is a refinement of the usual cycle index polynomial of $C_n$.

\begin{propos}[\cite{Timothy}]  The Walsh index series of $n$-cycles $C_n$ is given by
%\begin{eqnarray}
%W(C_n^\mathrm{or})=\frac{1}{n}\sum_{d|n}\varphi(d)a_d^{\frac{n}{d}}b_d^{\frac{n}{d}},
%\end{eqnarray}
%
\begin{eqnarray}
W_{C_n}({\bf a}; {\bf b}; {\bf c})=\frac{1}{2n}\sum_{d|n}\phi(d)a_d^{\frac{n}{d}}b_d^{\frac{n}{d}} + 
\frac{1}{2}\left\{ \begin{array}{ll} 
a_1a_2^{\frac{n-1}{2}}b_2^{\frac{n-1}{2}}c_1, & \mbox{\emph \,$n$ }\mathrm{odd}\\ 
\frac{1}{2}(a_2^{\frac{n}{2}}b_2^{\frac{n-2}{2}}c_1^2 + a_1^2a_2^\frac{n-2}{2}b_2^{\frac{n-2}{2}}), & \mbox{\emph \,$n$ }\mathrm{even} 
\end{array}
\right.
\label{form:cyc}
\end{eqnarray}
where $\phi$ is the Euler $\phi$-function.
\end{propos}

By combining the counting approaches of (\ref{form:bcp}) and (\ref{form:cyc}), we can compute the extended Walsh series for $C_n^\mathrm{m}$, $n\ge 3$, as follows. Notice that we use the homogenous %(two-variable) 
matching polynomials $U_n(y,z)$ and $T_n(y,z)$ %for cycles %$C_i$ of size $n$, 
given in Proposition 6.  The proof is left to the reader.
%which can also be found in \cite{GLL2} (see \cite{Godsil} as well).

\begin{thm}  The extended Walsh series $W_{C_n}^\mathrm{m}$ is given, for $n\geq3$, by the formulas
%
%\begin{eqnarray}
%W_{C_3}^\mathrm{m}({\bf a}; {\bf b}; {\bf c};\beta;\gamma)
%&=& \frac{1}{6}\left(2a_{{3}}\beta_{{3}}+{a_{{1}}}^{3} ({\beta_{{1}}}^{3}+3\,b_{{1}}{\beta_{{1}}}^{2} )\right)
%+\frac{1}{2}\,a_{{1}}a_{{2}} \left( \beta_{{2}}\gamma_{{1}}+\beta_{{2}}c_{{1}} \right),
%\label{form:bcc1}
%\end{eqnarray}
%
%for $n=4$, by
%
%\begin{eqnarray}
%W_{C_4}^\mathrm{m}({\bf a}; {\bf b}; {\bf c};\beta;\gamma)
%&=&
%\frac{1}{8}\left(2a_{{4}}\beta_{{4}}+{a_{{2}}}^{2} ({\beta_{{2}}}^{2}+2\,b_{{2}}\beta_{{2}} ) +
%{a_{{1}}}^{4} ( {\beta_{{1}}}^{4}+4\,b_{{1}}{\beta_{{1}}}^{3}+2\,{b_{{1}}}^{2}{\beta_{{1}}}^{2})\right) \nonumber\\
%&&+ \ \frac{1}{4}\left({a_{{1}}}^{2}a_{{2}}{\beta_{{2}}}^{2}+{a_{{2}}}^{2} ( {\gamma_{{1}}}^{2} ( b_{{2}}+\beta_{{2}} ) +2\,c_{{1}}\gamma_{{1}}\beta_{{2}}+{c_{{1}}}^{2}\beta_{{2}} )\right),
%\label{form:bcc2}
%\end{eqnarray}
%
\begin{eqnarray}
W_{C_n}^\mathrm{m}%(a_1,a_2,\ldots,a_n ;b_1,b_2,\ldots,b_n ;c_1 ;\beta_1,\beta_2,\ldots,\beta_n ; \gamma_1)
({\bf a}; {\bf b}; {\bf c};\beta;\gamma)
&=&\frac{1}{2n}\sum_{d|n}\phi(\frac{n}{d})a_\frac{n}{d}^dT_d(b_\frac{n}{d},\beta_\frac{n}{d})\nonumber \\
&&+\ \frac{1}{2}%\left\{ \begin{array}{ll}
a_1a_2^{\frac{n-1}{2}}(\beta_2\gamma_1U_\frac{n-1}{2}(b_2,\beta_2) + \beta_2^2c_1U_\frac{n-3}{2}(b_2,\beta_2)), \label{form:bcc3} \end{eqnarray}
for $n$ odd, %$n\ge 5$, 
and 
\begin{eqnarray}
W_{C_n}^\mathrm{m}%(a_1,a_2,\ldots,a_n ;b_1,b_2,\ldots,b_n ;c_1 ;\beta_1,\beta_2,\ldots,\beta_n ; \gamma_1)
({\bf a}; {\bf b}; {\bf c};\beta;\gamma)
&=&\frac{1}{2n}\sum_{d|n}\phi(\frac{n}{d})a_\frac{n}{d}^dT_d(b_\frac{n}{d},\beta_\frac{n}{d})
+\frac{1}{4}[a_1^2a_2^{\frac{n-2}{2}}\beta_2^2U_\frac{n-2}{2}(b_2,\beta_2)\nonumber \\
&+& a_2^\frac{n}{2}(\gamma_1^2U_\frac{n}{2}(b_2,\beta_2)
 + 2c_1\gamma_1\beta_2U_\frac{n-2}{2}(b_2,\beta_2)
+ c_1^2\beta_2^2U_\frac{n-4}{2}(b_2,\beta_2)) ], \nonumber \\
% \mbox{\emph \,$n\ge 6$ }\mathrm{even}. %\end{array}\right.
\label{form:bcc4}
\end{eqnarray}
for $n$ even.%, $n\ge 6$.
\end{thm}

It is then possible to obtain the generating function $\widetilde{C}_n^\mathrm{m}(x,y,z)$ of unlabelled matched cycles, for $n\ge 3$, where the variables $y$ and $z$ are edge-counters of sort $Y$ (in the matching) and $Z$ (not in the matching), respectively: set $a_i:=x^i, b_i:=y^i, c_i:=y^i, \beta_i:=z^i, \gamma_i:=z^i$ in (\ref{form:bcc3}) -- (\ref{form:bcc4}) for all possible values of the index $i$. Then we have the following corollary.

\begin{corol} The tilde generating function $\widetilde{C}_n^\mathrm{m}(x,y,z)$ of unlabelled matched cycles, for $n\ge 3$, is given by the following formulas:
\begin{eqnarray}
\widetilde{C}_n^\mathrm{m}(x,y,z) &=&\frac{1}{2n}\sum_{d|n}\phi(\frac{n}{d})x^nT_d(y^\frac{n}{d},z^\frac{n}{d})\nonumber \\
&+& 
\frac{1}{2}x^n(z^3U_\frac{n-1}{2}(y^2,z^2) + z^4yU_\frac{n-3}{2}(y^2,z^2)), 
\end{eqnarray}
for $n$ odd and 
\begin{eqnarray}
\widetilde{C}_n^\mathrm{m}(x,y,z) 
&=&\frac{1}{2n}\sum_{d|n}\phi(\frac{n}{d})x^nT_d(y^\frac{n}{d},z^\frac{n}{d})
+
\frac{1}{4}[x^n(z^4U_\frac{n-2}{2}(y^2,z^2)\nonumber \\
& + & z^2U_\frac{n}{2}(y^2,z^2) + 2yz^3U_\frac{n-2}{2}(y^2,z^2) + y^2z^4U_\frac{n-4}{2}(y^2,z^2))],  
\label{form:bccxyz}
\end{eqnarray}
for $n$ even.
\end{corol} 
%
%Note the special case
%%
%\begin{equation}
%\widetilde{C}_3^\mathrm{m}(x,y,z) 
%%&= &
%=\left( {z}^{3}+y{z}^{2} \right) {x}^{3},
%\end{equation}
%%
%\begin{equation}
%\widetilde{C}_4^\mathrm{m}(x,y,z) 
%%&= &
%=\left( {z}^{4}+y{z}^{3}+{y}^{2}{z}^{2} \right) {x}^{4},
%\end{equation}
%

Finally, we can compute the Walsh index series $W_\H({\bf a}; {\bf b}; {\bf c})$ of the species $\H$ of toroidal crowns by using (\ref{form:WH}), with $W_C^\mathrm{m}= \sum_{n\geq3}W_{C_n}^\mathrm{m}$ and $W^+_{K_5\backslash e}$ and 
$W^-_{K_5\backslash e}$ given by (\ref{form:w+k5}) and (\ref{form:w-k5}). 
%
%\section{Numerical results}
\section{Numerical results}
We used \textit{Maple IX} software to do all the computations. 
First we give the first terms of the tilde generating function $\widetilde{\H}(x,y)$ of unlabelled toroidal crowns, by setting $a_k:=x^k, b_k:=y^k, c_k:=y^k$ into (\ref{form:WH}): %$W_\H({\bf a}; {\bf b}; {\bf c})$:

\vspace{3mm}
\noindent
$\widetilde{\H}(x,y)=
{y}^{19}{x}^{9}+{y}^{20}{x}^{10}+{y}^{27}{x}^{12}+{y}^{28}{x}
^{13}+{y}^{29}{x}^{14}+{y}^{30}{x}^{15}+{y}^{36}{x}^{16}
+{y}^{37}{x}^{17}+2\,{y}^{38}{x}^{18}\\
+{y}^{39}{x}^{19}+  ({y}^{40}+{y}^{45})  {x}^{20}+{y}^{46}{x}^{21}+2\,{y}^{47}{x}^{22}+
2\,{y}^{48}{x}^{23}+  ({y}^{49}+{y}^{54})  {x}^{24}+  ({y}^{50}+{y}^{55})  {x}^{25}+3\,{y}^{56}{x}^{26}\\
+3\,{y}^{57}{x}^{27}+  (3\,{y}^{58}+{y}^{63})  {x}^{28}+
  ({y}^{59}+{y}^{64})  {x}^{29}+  ({y}^{60}+3\,{y}^
{65})  {x}^{30}+4\,{y}^{66}{x}^{31}+  (4\,{y}^{67}+{y}^
{72})  {x}^{32}\\
+  (3\,{y}^{68}+{y}^{73})  {x}^{33}
+  ({y}^{69}+4\,{y}^{74})  {x}^{34}+  ({y}^{70}+5\,
{y}^{75})  {x}^{35}+  (8\,{y}^{76}+{y}^{81})  {x}^
{36}+  (5\,{y}^{77}+{y}^{82})  {x}^{37}\\
+  (4\,{y}^{78}+4\,{y}^{83})  {x}^{38}+  ({y}^{79}+7\,{y}^{84})
  {x}^{39}+  ({y}^{80}+10\,{y}^{85}+{y}^{90})  {x}
^{40}+  (10\,{y}^{86}+{y}^{91})  {x}^{41}\\
+  (7\,{y}^{87}+5\,{y}^{92})  {x}^{42}+  (4\,{y}^{88}+8\,{y}^{93})
  {x}^{43}+  ({y}^{89}+16\,{y}^{94}+{y}^{99})  {x}
^{44}+  ({y}^{90}+16\,{y}^{95}+{y}^{100})  {x}^{45}
\\+  (16\,{y}^{96}+5\,{y}^{101})  {x}^{46}+  (8\,{y}^{
97}+10\,{y}^{102})  {x}^{47}+  (5\,{y}^{98}+20\,{y}^{103}
+{y}^{108})  {x}^{48}\\
+  ({y}^{99}+26\,{y}^{104}+{y}^{109})  {x}^{49}+  ({y}^{100}+26\,{y}^{105}+6\,{y}^{110})
  {x}^{50}+  (20\,{y}^{106}+12\,{y}^{111})  {x}^{51}\\
+  (10\,{y}^{107}+29\,{y}^{112}+{y}^{117})  {x}^{52}+
  (5\,{y}^{108}+38\,{y}^{113}+{y}^{118})  {x}^{53}+
  ({y}^{109}+50\,{y}^{114}+6\,{y}^{119})  {x}^{54}\\
  +  ({y}^{110}+38\,{y}^{115}+14\,{y}^{120})  {x}^{55}+
  (29\,{y}^{116}+35\,{y}^{121}+{y}^{126})  {x}^{56}+
  (12\,{y}^{117}+57\,{y}^{122}+{y}^{127})  {x}^{57}\\
  +  (6\,{y}^{118}+76\,{y}^{123}+7\,{y}^{128})  {x}^{58}+
  ({y}^{119}+76\,{y}^{124}+16\,{y}^{129})  {x}^{59}\\
  +  ({y}^{120}+57\,{y}^{125}+47\,{y}^{130}+{y}^{135})  {x}^
{60}+  (35\,{y}^{126}+79\,{y}^{131}+{y}^{136})  {x}^{61}\\
+  (14\,{y}^{127}+126\,{y}^{132}+7\,{y}^{137})  {x}^{62}+
  (6\,{y}^{128}+133\,{y}^{133}+19\,{y}^{138})  {x}^{63}\\
  +  ({y}^{129}+126\,{y}^{134}+56\,{y}^{139}+{y}^{144})  {x}^{64}+ ...
%  ({y}^{130}+79\,{y}^{135}+111\,{y}^{140}+{y}^{145}) {x}^{65})
$

\vspace{3mm}

Recall that $\T_C=K_5 + M + M^* + \H$. It suffices to add to the above series the terms 
$$x^5y^{10} + (y^{19} +y^{18})x^8,$$ 
corresponding to $K_5$, $M$ and $M^*$, to cover all toroidal cores with up to 64 vertices. By setting $y=1$, we obtain the numbers $t_C(n)$ of unlabelled toroidal cores with $n$ vertices, presented in Table~1. 

%
%Table 1
\begin{table}[!h] \label{table:Tc}
\centerline {\scriptsize
	\begin{tabular}[t]{|| r | r | r | r | r | r | r | r | r | r | r | r | r | r | r | r | r | r | r | r | r ||}
	\hline
	$n$ & 5 & 6 & 7 & 8 & 9 & 10 & 11 & 12 & 13 & 14 & 15 & 16 & 17 & 18 & 19 & 20 & 21 & 22 & 23 & 24\\
	\hline
	$t_C(n)$ & 1 & 0 & 0 & 2 & 1 & 1 & 0 & 1 & 1 & 1 & 1 & 1 & 1 & 2 & 1 & 2 & 1 & 2 & 2 & 2\\
	\hline
	\hline
	$n$ & 25 & 26 & 27 & 28 & 29 & 30 & 31 & 32 & 33 & 34 & 35 & 36 & 37 & 38 & 39 & 40 & 41 & 42 & 43 & 44\\
	\hline
	$t_C(n)$ & 2 & 3 & 3 & 4 & 2 & 4 & 4 & 5 & 4 & 5 & 6 & 9 & 6 & 8 & 8 & 12 & 11 & 12 & 12 & 18\\
	\hline
	\hline
	$n$ & 45 & 46 & 47 & 48 & 49 & 50 & 51 & 52 & 53 & 54 & 55 & 56 & 57 & 58 & 59 & 60 & 61 & 62 & 63 & 64\\
	\hline
	$t_C(n)$ & 18 & 21 & 18 & 26 & 28 & 33 & 32 & 40 & 44 & 57 & 53 & 65 & 70 & 89 & 93 & 106 & 115 & 147 & 158 & 184\\
	\hline
	\end{tabular} }
	\caption{The number $t_C(n)$ of unlabelled toroidal cores (having $n$ vertices).}
\end{table}

To count unlabelled toroidal graphs with no $K_{3,3}$'s, it is necessary to know the generating functions for unlabelled strongly planar networks $\N_P$ and $\tau$-symmetric networks $\N_{P,\tau}$. It is still an open problem to give these power series in general. %without generating any maps. 
We have enumerated planar networks with up to six (four internal) vertices by hand. The corresponding generating functions are:
\begin{eqnarray}
\widetilde{\N}_P(x,y) & = & y+(1+y)[ y^2x+(y^3+3y^4+y^5)x^2+(y^4+8y^5+15y^6+9y^7+3y^8)x^3\nonumber \\
 & & +(y^5 + 16y^6 + 66y^7+ 112y^8 + 97y^9 + 47y^{10} + 9y^{11})x^4+\ldots] ,
\label{form:NP}
\end{eqnarray}
and 
\begin{eqnarray}
\widetilde{\N}_{P,\tau}(x,y) & = & y+(1+y)[y^2x+(y^3+y^4+y^5)x^2+(y^4+2y^5+3y^6+3y^7+y^8)x^3\nonumber \\
 & & +(y^5 + 4y^6 + 8y^7+ 12y^8 + 13y^9 + 7y^{10} + 3y^{11})x^4+\ldots].
\label{form:NPtau}
\end{eqnarray}

Tables~2 and~3 present the corresponding results for the classes $\P$ and $\T-\P$ of non-planar 2-connected $K_{3,3}$-free projective-planar and non-projective-planar toroidal graphs. Table~2 gives the number of graphs in $\P$, and Table~3 gives the number of graphs in $(\T-\P)$. Recall that $\P = K_5\uparrow\N_P$ and $\T = \T_C\uparrow\N_P$. Therefore $(\T-\P) = (\T_C-K_5)\uparrow\N_P$. 

%
%Table 2
\begin{table}[!h] \label{table:fnm}
\centerline {\footnotesize
	\begin{tabular}[t]{|| r | r | r || r | r | r || r | r | r ||}
	\hline
	$n$ & $m$ & $f_{n,m}$ & $n$ & $m$ & $f_{n,m}$ & $n$ & $m$ & $f_{n,m}$\\
	\hline \hline
	5 & 10 & 1 &  8 & 13 & 7 &  9 & 14 & 17\\
	\hline
	6 & 11 & 1 & 8 & 14 & 21 &  9 & 15 & 76\\
	\hline
	6 & 12 & 1 & 8 & 15 & 34 &  9 & 16 & 197\\
	\hline
	7 & 12 & 3 & 8 & 16 & 28 &  9 & 17 & 272\\
	\hline
	7 & 13 & 5 & 8 & 17 & 10 &  9 & 18 & 234\\
	\hline
	7 & 14 & 5 & 8 & 18 & 2 &  9 & 19 & 120\\
	\hline
	7 & 15 & 1 & \multicolumn{3}{c ||}{ } &  9 & 20 & 40\\
	\cline{1-3}\cline{7-9}
	\multicolumn{6}{c ||}{ } &  9 & 21 & 6\\
	\cline{7-9}
\end{tabular} }
\caption{The number of unlabelled 2-connected $K_{3,3}$-free non-planar projective-planar graphs %with no $K_{3,3}$'s 
with $n$ vertices and $m$ edges.}
\end{table}
%

%
%Table 3
\begin{table}[!h] \label{table:fnm}
\centerline {\footnotesize
	\begin{tabular}[t]{|| r | r | r || r | r | r || r | r | r ||}
	\hline
	$n$ & $m$ & $f_{n,m}$ & $n$ & $m$ & $f_{n,m}$ & $n$ & $m$ & $f_{n,m}$\\
	\hline \hline
	8 & 18 & 1 &  11 & 21 & 67 & 12 & 22 & 277\\
	\hline
	8 & 19 & 1 & 11 & 22 & 245 & 12 & 23 & 1361\\
	\hline
	9 & 19 & 3 & 11 & 23 & 419 & 12 & 24 & 3274\\
	\hline
	9 & 20 & 5 & 11 & 24 & 396 & 12 & 25 & 4598\\
	\hline
	9 & 21 & 3 & 11 & 25 & 204 & 12 & 26 & 4061\\
	\hline
	10 & 20 & 17 & 11 & 26 & 50 & 12 & 27 & 2295\\
	\hline
	10 & 21 & 39 & 11 & 27 & 7 & 12 & 28 & 823\\
	\hline
	10 & 22 & 44 & \multicolumn{3}{c ||}{ } & 12 & 29 & 195\\
	\cline{1-3}\cline{7-9}
	10 & 23 & 24 & \multicolumn{3}{c ||}{ } & 12 & 30 & 21\\
	\cline{1-3}\cline{7-9}
	10 & 24 & 3 & \multicolumn{6}{c}{ }\\
	\cline{1-3}
	\end{tabular} }
\caption{The number of unlabelled 2-connected $K_{3,3}$-free non-projective-planar toroidal graphs %with no $K_{3,3}$'s 
with $n$ vertices and $m$ edges.}
\end{table}

Timothy Walsh has extended the series (\ref{form:NP}) and (\ref{form:NPtau}) to eight internal vertices with the help of Brendon McKay, author of the software ``plantri'' (available on the web, see \cite{BMplantri}) who supplied a list of all the 3-connected planar (embedded) graphs and their automorphisms and by using techniques of \cite{Timothy}. 
Consequently, Tables 2 and 3 have been extended to $n$ = 13 and 16, respectively, and are available from the authors. Further extensions are possible.

Since vertices of degree two are uninteresting %meaningless 
for graph embeddability questions, we plan to count the corresponding homeomorphically irreducible graphs in a future work.% with no vertices of degree two.
%
%\section{Proofs}
%
\section{Proofs}
Our principal goal in this section is to prove Theorems\,3 and\,4. We also provide proofs of the other results and formulas of Section\,3. 
The proof given by Walsh of his Proposition 2.2 (see \cite{Timothy}), which implies Theorem 3, is rather sketchy. It is based on a result on the cycle index polynomial of wreath products used by Robinson to prove his composition theorem for graphs (see \cite{Robinson}). Walsh gives two examples to show how the approach of Robinson can be adapted to prove his proposition. 

Here we provide complete and detailed proofs based on the methods of species theory (see for example, \cite{BLL}, section 4.3). 
To achieve this, we reformulate the Walsh index series $W_\G$ and the (tilde) generating function $%\widetilde{\G}%
\G^\sim(x,y)$ in terms of labelled enumeration of associated weighted species. Following an idea of Joyal \cite{Joyal}, we introduce the auxiliary weighted species $\G^\aut=\G^\aut_w$. For any finite set $U$ (of vertices), $\G^\aut[U]$ is defined as the set of graphs in $\G[U]$ equipped with an automorphism $\sigma$, i.e. 
\[
\G^\aut[U]=\{ (G,\sigma)\ |\ G\in \G[U], \sigma\in S[U] : \sigma\cdot G=G \},
\]
where $S[U]$ is the set of all permutations of $U$. The relabelling rule of $\G^\aut$-structures along a bijection $\beta : U\ \tilde{\longrightarrow}\ U^\prime$ is defined as follows:
\[
\beta \cdot (G,\sigma) = (\beta \cdot G,\ \beta {\circ} \sigma {\circ} \beta^{-1}),
\]
where $\beta \cdot G$ is the graph obtained from $G$ by relabelling along $\beta$ and the composition $\circ$ is taken from right to left. It is easy to verify that $\G^\aut_w$ is a well-defined weighted species, where the weight function $w(G,\sigma)$ is the cycle index monomial defined by (\ref{form:weightGsigma}). Recall that $|\G^\aut[n]|_w$ denotes the total weight of $\G^\aut_w$-structures over the vertex set $[n]:=\{1,2,\ldots,n\}$, i.e. 
\begin{eqnarray*}
|\G^\aut[n]|_w=\sum_{(G,\sigma)\in \G^\aut_w[n]} w(G,\sigma).
\label{form:WabcAut}
\end{eqnarray*}

We will also use the weight function $w_0$ defined as $w_0(G,\sigma)=y^m$, where $m$ is the number of edges in $G$. 
As the following two propositions show, the Walsh index series $W_\G({\bf a}; {\bf b}; {\bf c})$ and the  generating function $\widetilde{\G}(x,y)$ appear as special cases of the usual exponential generating function of weighted species 
\begin{eqnarray*}
\G_w^\aut(x) = \sum_{n\ge 0}|\G^\aut[n]|_w\frac{x^n}{n!},
\label{form:WabcAut}
\end{eqnarray*}
with the weight functions $w$ as above, and $w=w_0$, respectively.
%Two special cases are given in Propositions\,9 and\,10.

\begin{propos}
Using the exponential generating function of labelled $\G^\aut_w$-structures, we have
\begin{eqnarray}
\G^\aut_w(x)|_{x=1} = W_\G({\bf a}; {\bf b}; {\bf c}).
\end{eqnarray}
\end{propos}
\textit{Proof.} The number of distinct graphs on $[n]$ obtained by relabelling a given graph $G$ with $n$ vertices is given by $\frac{n!}{|\mathrm{Aut}(G)|}$. Therefore formula (\ref{form:Wabc}) gives:
\begin{eqnarray*}
W_\G({\bf a}; {\bf b}; {\bf c}) & = & \sum_{n\ge 0} \sum_{G\in \mathrm{Typ}(\G_n)} \frac{1}{|\mathrm{Aut}(G)|} \sum_{\sigma\in \mathrm{Aut}(G)} w(G,\sigma)\nonumber \\
 & = & \sum_{n\ge 0} \frac{1}{n!} \sum_{G\in \mathrm{Typ}(\G_n)} \frac{n!}{|\mathrm{Aut}(G)|} \sum_{\sigma\in \mathrm{Aut}(G)} w(G,\sigma)\nonumber \\
 & = & \sum_{n\ge 0} \frac{1}{n!} \sum_{G\in \G[n]} \sum_{\sigma\in \mathrm{Aut}(G)} w(G,\sigma)\nonumber \\
 & = & \sum_{n\ge 0} \frac{1}{n!} \sum_{(G,\sigma)\in \G^\aut[n]} w(G,\sigma)\nonumber \\
 & = & \sum_{n\ge 0} \frac{1}{n!}|\G^\aut[n]|_w = \G^\aut_w(x)|_{x=1}.
\end{eqnarray*}
\hfill\rule{2mm}{2mm}

\begin{propos}
Using the exponential generating function of labelled $\G^\aut_{w_0}$-structures, we have
\begin{eqnarray}
\G^\aut_{w_0}(x) = \widetilde{\G}(x,y).
\end{eqnarray}
\end{propos}
\textit{Proof.} It is easy to see, using Burnside's Lemma (alias Cauchy-Frobenius formula, for example, see \cite{BLL}), that
\begin{eqnarray*}
\G^\aut_{w_0}(x) & = & \sum_{n\ge 0} \sum_{(G,\sigma)\in \G^\aut[n]} w_0(G,\sigma) \frac{x^n}{n!} \nonumber \\
& = & \sum_{n\ge 0} \frac{1}{n!} \sum_{\sigma\in S_n} \sum_{G\in \mathrm{Fix}\,\G[\sigma]} x^ny^m \nonumber \\
& = & \sum_{n\ge 0} \sum_{G\in \mathrm{Typ}\,\G[n]} x^ny^m \nonumber \\
& = & \widetilde{\G}(x,y).
\end{eqnarray*}
Here the notation $\mathrm{Fix}\,\G[\sigma]=\{G\in \G[n]\,|\,\sigma\cdot G=G \}$ has been used.
\hfill\rule{2mm}{2mm}

\bigbreak
\noindent\textbf{Proof of Proposition 1.} The left-hand side of (\ref{form:7}) is the mixed exponential generating function that counts labelled graphs in $\G$. After the substitution in the right-hand side of (\ref{form:7}), following the proof of Proposition\,10, the only surviving permutation is $\sigma=\mathrm{id}_{[n]}$, and the weight function becomes $w(G,\sigma)=x^ny^m$, where $n$ is the number of vertices and $m$ is the number of edges in $G$. Then we have
\begin{eqnarray*}
W_\G(x,0,0,\ldots; y,0,0,\ldots; 0,0,0,\ldots) & = & \sum_{n\ge 0} \frac{1}{n!} \sum_{G\in \G[n]} \sum_{\sigma\in \mathrm{Aut}(G)} w(G,\sigma)\nonumber \\ 
 & = & \sum_{n\ge 0} \frac{1}{n!} \sum_{G\in \G[n]} x^ny^m \nonumber \\ 
 & = & \G(x,y).
\end{eqnarray*}

To prove (\ref{form:8}), note that 
$$\bar{w}(G,\sigma):=w(G,\sigma)|_{a_i=x^i,\, b_i=y^i,\, c_i=y^i}=x^ny^m$$ 
and, using Propositions\,10 and\,11, we have
$$W_\G(x,x^2,x^3,\ldots; y,y^2,y^3,\ldots; y,y^2,y^3,\ldots)  = \G^\aut_{w}(x)|_{x=1, w=\bar{w}} = \G^\aut_{w_0}(x) = 
\widetilde{\G}(x,y).$$
\hfill\rule{2mm}{2mm}

A similar approach can be used for the Walsh index series $W_\N^+$ and $W_\N^-$ of a given species of $2$-pole networks $\N$. Denote by 
$$\N^+[U]=\{ (N,\sigma)\ |\ N\in \N[U],\ \sigma\in \mathrm{Aut}^+(N) \}$$ 
and 
$$\N^-[U]=\{ (N,\sigma)\ |\ N\in \N[U],\ \sigma\in \mathrm{Aut}^-(N) \},$$
where $\Aut^+(N)$ and $\Aut^-(N)$ are defined by (\ref{form:A+}) and (\ref{form:A-}), respectively.
Then, using the weight functions given by (\ref{form:weightN+}) and (\ref{form:weightN-}),
$\N^+_w$ and $\N^-_w$ are weighted species whose labelled enumerations yield by specialization the series $W^+_\N$, $W^-_\N$, $\widetilde{\N}(x,y)$ and $\widetilde{\N}_\tau(x,y)$.

\begin{propos}
For a species of networks $\N$, the Walsh index series $W^+_\N$ and $W^-_\N$ can be expressed by the formulas
\begin{eqnarray}
\N^+_w(x) |_{x=1} = W^+_\N({\bf a}; {\bf b}; {\bf c})
\label{form:34}
\end{eqnarray}
and
\begin{eqnarray}
\N^-_w(x) |_{x=1} = W^-_\N({\bf a}; {\bf b}; {\bf c})
\label{form:35}
\end{eqnarray}
respectively.
\end{propos}
\textit{Proof.} Notice that for (\ref{form:34}) we have $\Aut^+(N)=\Aut(N)$, and, when $N\in \N_\tau$, then for (\ref{form:35}) we have $|\Aut^-(N)|=|\Aut^+(N)|$, since given any $\varphi\in \Aut^-(N)$, the application $\sigma \mapsto \sigma\circ\varphi$ defines a bijection from $\Aut^+(N)$ to $\Aut^-(N)$. Then, by using Definition\,2, the proof is similar to the proof of Proposition\,10.
\hfill\rule{2mm}{2mm}

\medbreak
Defining the weight function $w_0$ as $w_0(N,\sigma)=y^m$ for a network $N$ having $m$ edges, we obtain the following identities.
\begin{propos}
Using the usual (exponential) generating function of labelled $\N^+_{w_0}$-structures and $\N^-_{w_0}$-structures, we have
\begin{eqnarray}
\N^+_{w_0}(x) = \widetilde{\N}(x,y)
\end{eqnarray}
and
\begin{eqnarray}
\N^-_{w_0}(x) = \widetilde{\N}_\tau(x,y)
\end{eqnarray}
respectively.
\end{propos}
\textit{Proof.} Using the remarks of the proof of Proposition\,12, the proof here is analogous to that of Proposition\,11.
\hfill\rule{2mm}{2mm}

\bigbreak
\noindent\textbf{Proof of Proposition 2.} Notice that 
$$w(N,\sigma)|_{a_i=x^i,\, b_i=y^i,\, c_i=y^i} = x^ny^m,$$
i.e.
$$\frac{w(\hat{N},\sigma^+)}{a^2_1}|_{a_i=x^i,\, b_i=y^i,\, c_i=y^i} = x^ny^m$$
and 
$$\frac{w(\hat{N},\sigma^-)}{a_2}|_{a_i=x^i,\, b_i=y^i,\, c_i=y^i} = x^ny^m,$$
where $n$ is the number of internal vertices and $m$ is the number of edges of $N$.
Then the proof is similar to that of Proposition\,1, by using Definition\,2. 
\hfill\rule{2mm}{2mm}

\bigbreak
\noindent\textbf{Proof of Proposition 3.} 
In order to prove (\ref{eq:WB01+}), it is sufficient to establish the combinatorial equality of weighted species
\begin{equation}\label{eq:partialw}
2\frac{\partial}{\partial b_1}\B^\aut_w = a_1^2X^2 \cdot \B^+_{01,w}.
\end{equation}
Now, to obtain a $2\frac{\partial}{\partial b_1}\B^\aut_w$-structure on a set $U$, we start with a $\B^\aut_w$-structure $(B,\sigma)$ on $U$; the operator $2\frac{\partial}{\partial b_1}$ is then interpreted as selecting, orienting and putting aside a cylindrical edge cycle of $\sigma$ of length 1, $\vec{e}=(u,v)$, where $\{u,v\}\subseteq U$; it is natural to relabel the vertices $u$ and $v$ by $0$ and $1$, respectively, in the graph $B$.  This is equivalent to first selecting and ordering the two vertices $u$ and $v$ from $U$ and defining $\sigma$ to be the identity on $\{u,v\}$, which yields an $a_1^2X^2$-structure, and then selecting a $\B^+_{01,w}$-structure $(B_0,\sigma_0)$ on the complementary set $U_0=U\backslash \{u,v\}$, giving the second factor on the right hand side of (\ref{eq:partialw}). The reverse construction consists simply of setting $U=U_0\cup \{u,v\}$, $B=B_0$, identifying $u$ with $0$ and $v$ with $1$, and $\sigma=\sigma_0\cup\mathrm{id}_{\{u,v\}}$. Taking generating functions and setting $x=1$ yields $$2 \frac{\partial}{\partial b_1}W_{\B}({\bf a}; {\bf b}; {\bf c})=a_1^2W_{\B_{0,1}}^+({\bf a}; {\bf b}; {\bf c})$$ which is equivalent to (\ref{eq:WB01+}). 
Also, multiplication of $\B^+_{01,w}$ by $b_1$ corresponds to reinserting the edge $e$  in $B$ so that we have an isomorphism 
\begin{equation}\label{form:57}
\N_{\B,w}^+ = (1+b_1) \B^+_{01,w}- \nn
\end{equation}
which yields (\ref{form:20}).

Notice that the same species isomorphisms, with the weight $w_0$, would give (\ref{eq:B01xy}) and (\ref{eq:NBxy}).  
The proof of (\ref{form:19}) and (\ref{form:21}) is based on the species isomorphism
\begin{equation}\label{form:59}
2\frac{\partial}{\partial c_1}\B^\aut_w = a_2X^2 \cdot \B^-_{01,w}
\end{equation}
which is established in a similar manner. 
\hfill\rule{2mm}{2mm}

\bigbreak
For a species of graphs $\G$, we introduce the following notation which will be used in the proofs of Theorems\,3 and\,4. Given a permutation $\sigma\in S_n$, we denote by $\underline{n}=(n_1,n_2,n_3\ldots)$ the cycle type of $\sigma$. 
Notice that $n_1+2n_2+3n_3+\ldots = n < \infty$, which we write simply as $\underline{n} < \infty$, and that the number of permutations of a given cycle type $\underline{n}$ is equal to $\frac{n!}{1^{n_1}n_1!2^{n_2}n_2!\ldots}$. 
If $G\in\mathrm{Fix}\,G[\sigma]$, denote by $k_i=k_i(G,\sigma)$ and $m_i=m_i(G,\sigma)$ the number of cylindrical and M\"obius edge cycles, respectively, of length $i$ induced by $\sigma$. Notice that the $k_i$'s and $m_i$'s are completely determined by the cycle type of the automorphism $\sigma$ (for example, see \cite{Timothy, RW1, RW2, BLL}). Therefore we can write $k_i=k_i(G,\underline{n})$ and $m_i=m_i(G,\underline{n})$. We also set 
$$\mathrm{Fix}\,\G[\underline{n}] = \mathrm{Fix}\,G[\sigma],
$$
where $\sigma$ is some fixed permutation of type $\underline{n}$. Then we have the following alternate expression for $W_\G({\bf a}; {\bf b}; {\bf c})$:
\begin{eqnarray}
W_\G({\bf a}; {\bf b}; {\bf c}) & = & \sum_{n\ge 0} \frac{1}{n!} \sum_{(G,\sigma)\in \G^\aut[n]} w(G,\sigma)\nonumber \\
 & = & \sum_{n\ge 0} \frac{1}{n!} \sum_{\sigma\in S_n} \sum_{G\in \mathrm{Fix}\,G[\sigma]} a_1^{n_1}a_2^{n_2}\ldots b_1^{k_1}b_2^{k_2}\ldots c_1^{m_1}c_2^{m_2}\ldots \nonumber \\
 & = & \sum_{\underline{n}<\infty}  \frac{1}{1^{n_1}n_1!2^{n_2}n_2!\ldots} \sum_{G\in \mathrm{Fix}\,\G[\underline{n}]} a_1^{n_1}a_2^{n_2}\ldots b_1^{k_1}b_2^{k_2}\ldots c_1^{m_1}c_2^{m_2}\ldots. \label{eq:WGter}%\nonumber \\
\end{eqnarray}
\noindent\textbf{Example.} As an illustration, let us consider the class $K_n$ of complete graphs, $n\ge 1$. Its Walsh index series is given by (see \cite{RW1, RW2, Timothy}, and also \cite{BLL}).
%
%a of the cycle type of edges induced and completely defined by the cycle type of a permutation of the graph vertices 
% Suppose the the permutation $\sigma$ of the complete graph $K_n$ vertices has the cycle type $\underline{n}=(n_1,n_2,n_3\ldots)$, $n_1+2n_2+3n_3+\ldots = n < \infty$. Then the Walsh index series of $K_n$ is given by
%
\begin{eqnarray}
W_{K_n}({\bf a}; {\bf b}; {\bf c}) &=& \sum_{\underline{n}} 
\prod_i \frac{a_i^{n_i}}{i^{n_i}n_i!} 
\prod_{i<j}b_{[i,j]}^{(i,j)n_in_j} 
\prod_{i}b_i^{i{n_i\choose2}+\lfloor\frac{(i-1)}{2}\rfloor n_i} 
\prod_{i}c_i^{n_{2i}},
\label{form:Wkn}
\end{eqnarray}
where the sum is taken over all sequences $\underline{n}=(n_1,n_2,n_3\ldots)$ such that $n_1+2n_2+3n_3+\ldots = n$ and where
$[i,j]$ denotes the least common multiple and $(i,j)$, the greatest common divisor of $i$ and $j$. 
Notice that the substitution $a_i:=a_i, b_i:=1+b_i, c_i:=1+c_i$ into (\ref{form:Wkn}) gives the Walsh index series $W_{\G_n}$ for the species of all graphs on $n$ vertices.
In other words, we have 
\begin{eqnarray} \label{form:WGn}
W_{\G_{\mathrm{a}}}({\bf a}; {\bf b}; {\bf c})&=&W_{K}({\bf a}; {\bf 1+b}; {\bf 1+c}),
\end{eqnarray}
where $\G_{\mathrm{a}}$ denotes the species of all graphs and $K$, the species of complete graphs.
As observed in \cite{Timothy}, this can be seen as an application of Theorem 3 since in fact we have 
\begin{eqnarray} \label{form:Ga=Kflecheun}
\G_{\mathrm{a}}&=&K\uparrow\N_{K_2},
\end{eqnarray}
where $\N_{K_2}$, denotes the class of trivial networks consisting of either two isolated poles or two poles joined by an edge.%, the composition being canonical.
%See also \cite{RW1, RW2, Timothy}.

\smallskip%\bigbreak
Given networks $N_1\in \N[U]$ and $N_2\in \N[V]$, 
%a mapping $\varphi : N_1 \longrightarrow N_2$ which is 
a bijection between the underlying sets (of internal vertices) $\varphi : U \longrightarrow V$ is a \textit{network isomorphism} if 
$$\varphi^+ = \varphi \cup \mathrm{id}_{\{0,1\}}: \hat{N}_1 \longrightarrow \hat{N}_2$$ 
is a graph isomorphism, where $\hat{N}$ denotes the graph on $U\cup\{0,1\}$ associated to the network $N\in\N[U]$.
Similarly, a bijection $\varphi : U \longrightarrow V$ is a \textit{network anti-isomorphism} if 
$$\varphi^- = \varphi \cup \tau: \hat{N}_1 \longrightarrow \hat{N}_2$$ 
is a graph isomorphism, where $\tau$ denotes the tranposition $(0,1)$.
%is an operator that maps pole $0$ of $N_1$ to the pole $1$ of $N_2$ and pole $1$ of $N_1$ to the pole $0$ of $N_2$.

\smallskip
A \textit{cylindrical $m$-wreath of networks} %of length $m$, 
%denoted by $\K_m(\N_{w_0})$, 
is defined as an oriented cycle $c_m$ of length $m$ of network isomorphisms 
\begin{equation}\label{eq:cm}
c_m : N_1 \stackrel{\varphi_1}{\longrightarrow} N_2 \stackrel{\varphi_2}{\longrightarrow}\ldots\stackrel{\varphi_{m-2}}{\longrightarrow} N_{m-1} \stackrel{\varphi_{m-1}}{\longrightarrow} N_m \stackrel{\varphi_m}{\longrightarrow} N_1,
\end{equation} 
where the networks $N_i$ are assumed to have disjoint %underlying 
sets of vertices. 
%and the respective poles are assumed to be distinct.
%

Let $\K_m(\N)$ denote the species of cylindrical $m$-wreaths of networks.
%of length $m$.
We associate the weights $w(c_m)$ and $w_0(c_m)$ as follows. Define the graph $\widehat{c_m}$ to be the %sum (disjoint union) 
disjoint union of all the graphs $\widehat{N_i}$ associated to the networks $N_i$ appearing in $c_m$. 
Also define $\sigma^+$ to be the graph automorphism of $\widehat{c_m}$ given by 
$$
\sigma^+ = \cup_{i=1}^m\,\varphi_i^+.
$$
Then we set $w(c_m):= w(\widehat{c_m},\sigma^+)/{a_m^2}$, where the weight $w(\widehat{c_m},\sigma^+)$ %$w(G,\sigma)$ 
is given by (\ref{form:weightGsigma}), and $w_0(c_m):=y^{|E(\widehat{c_m})|}$.

A \textit{rooted} cylindrical $m$-wreath of networks %of length $m$ 
is a cylindrical $m$-wreath of networks $c_m$ where one network is distinguished from the others.
In fact the description (\ref{eq:cm}) of $c_m$ includes a rooting at $N_1$. In the unrooted case, all the possible rootings are considered equivalent. Let $\K_m^{\bullet}(\N)_{w}$ and $\K_m^{\bullet}(\N)_{w_0}$ denote the corresponding weighted species of rooted cylindrical $m$-wreaths of networks.
We introduce the weights $(w)_m$ and $w_0^m$ on $\N^+$-structures $(N,\sigma)$  by 
\begin{equation} \label{eq:wmw0m}
(w)_m(N,\sigma)= w(N,\sigma)|_{a_i:=a_{mi},a_i:=a_{mi},a_i:=a_{mi}} \ \ 
\mathrm{and} \ \  w_0^m(N,\sigma)=w_0(N,\sigma)^m.
\end{equation}
We then have, using the plethystic notation (\ref{eq:W+Nk}),
\begin{equation} \label{eq:N+wm0m}
\N^+_{(w)_m}(x)|_{x=1} = (W^+_\N)_m%({\bf a}; {\bf b}; {\bf c})
 \ \ 
\mathrm{and} \ \  \N^+_{w_0^m}(x)= \widetilde{\N}(x,y^m).
\end{equation}
\begin{propos} We have the weighted species isomorphisms 
\begin{equation}\label{eq:KmNbullet}
\K_m^{\bullet}(\N)_{w} =%\cong 
\N_{(w)_m}^+(X^m), \ \ \ \K_m^{\bullet}(\N)_{w_0} =%\cong 
\N_{w_0^m}^+(X^m)
\end{equation}
and the series equalities
\begin{equation}\label{eq:KmNbulletx}
\K_m^{\bullet}(\N)_{w}(x) = \N_{(w)_m}^+(x^m), \ \ \ 
\K_m^{\bullet}(\N)_{w_0}(x) = \widetilde{\N}(x^m,y^m).
\end{equation}
\end{propos}
\noindent\textit{Proof}.  
Given a rooted cylindrical $m$-wreath of networks $c_m$ in $\K_m^{\bullet}(\N)$, of the form (\ref{eq:cm}), note that the composite 
$\varphi_0 = \varphi_m \circ \varphi_{m-1} \circ \ldots \varphi_2 \circ \varphi_1$ is an automorphism of $N_1$, and we obtain a $\N^+$-structure $(N_1,\varphi_0)$. Moreover the sequence of network isomorphisms $(\varphi_1,\ldots,\varphi_{m-1})$ can be encoded in a set of lists of length $m$ $(u_1, u_2, \dots, u_m)$, where $u_1$ runs over the underlying set of $N_1$ and $u_{i+1}= \varphi_i(u_i)$, $i=1\ldots m-1$, and we can consider the $\N^+$-structure $(N_1,\varphi_0)$ to ``live" on this set of lists. In other words, what we have obtained is an $\N^+(X^m)$-structure.  Since the isomorphism $\varphi_m$ can be recovered from $\varphi_0$ and the other isomorphisms by the rule $\varphi_m = \varphi_0 \circ \varphi_1^{-1} \circ \varphi_2^{-1} \circ \ldots \circ \varphi_{m-1}^{-1}$, this correspondence is bijective. Finally we can see that the vertex and edge cycle structure of $(N_1,\varphi_0)$ is the same as that of $(\widehat{c_m},\sigma^+)$, except for the fact that all cycle lengths are multiplied by $m$ in $(\widehat{c_m},\sigma^+)$, and that the number of edges of $N_1$ is multiplied by $m$ in $\widehat{c_m}$.
This proves the two combinatorial identities of (\ref{eq:KmNbullet}). The formulas of  (\ref{eq:KmNbulletx}) then follow by taking generating functions.
\hfill\rule{2mm}{2mm}

\medskip%\bigbreak
We now define a \textit{rooted M\"obius $m$-wreath of networks} $c^a_m$ to be a sequence of network isomorphisms 
$N_1 \stackrel{\varphi_1}{\longrightarrow} N_2 \stackrel{\varphi_2}{\longrightarrow}\ldots\stackrel{\varphi_{m-2}}{\longrightarrow} N_{m-1} \stackrel{\varphi_{m-1}}{\longrightarrow} N_m$ of length $m-1$ 
followed by a network anti-isomorphism 
$\varphi_m : N_m \longrightarrow N_1$. Notice that 
$\varphi_0 = \varphi_m \circ \varphi_{m-1} \circ \ldots \varphi_2 \circ \varphi_1$ is an anti-automorphism of $N_1$, and that
$\varphi_m = \varphi_0 \circ \varphi_1^{-1} \circ \varphi_2^{-1} \circ \ldots \circ \varphi_{m-1}^{-1}$.
%denoted by  %or $\K_m^{\mu\bullet}(\N_{w_0})$
%of length $m$.
%We associate the weights $w(c_m)$ and $w_0(c_m)$ as follows. 
Define the graph $\widehat{c^a_m}$ to be the disjoint union of all the graphs $\widehat{N_i}$ associated to the networks $N_i$. %appearing in $c^a_m$. 
Also define $\sigma^-$ to be the graph automorphism of $\widehat{c^a_m}$ given by 
$$
\sigma^- = \cup_{i=1}^{m-1}\,\varphi_i^+\cup\varphi_m^-.
$$
Then we set $w(c^a_m):= w(\widehat{c^a_m},\sigma^-)/{a_m^2}$ %where the weight $w(\widehat{c_m},\sigma^+)$ %$w(G,\sigma)$ is given by (\ref{form:weightGsigma}), 
and $w_0(c^a_m):=y^{|E(\widehat{c^a_m})|}$.
Let $\K_m^{\bullet\mu}(\N)_w$ and $\K_m^{\bullet\mu}(\N)_{w_0}$ denote the corresponding species of weighted rooted M\"obius $m$-wreaths of networks of $\N$.
\begin{propos} We have the weighted species isomorphisms 
\begin{equation}\label{eq:KmmuNbullet}
\K_m^{\bullet\mu}(\N)_{w} =%\cong 
\N_{(w)_m}^-(X^m), \ \ \ \K_m^{\bullet\mu}(\N)_{w_0} =%\cong 
\N_{w_0^m}^-(X^m)
\end{equation}
and the series equalities
\begin{equation}\label{eq:KmmuNbulletx}
\K_m^{\bullet\mu}(\N)_{w}(x) = \N_{(w)_m}^-(x^m), \ \ \ 
\K_m^{\bullet\mu}(\N)_{w_0}(x) = \widetilde{\N}_\tau(x^m,y^m).
\end{equation}
\end{propos}
\noindent\textit{Proof}. The proof is similar to that of Proposition 14. Details are left to the reader.
\hfill\rule{2mm}{2mm}

\smallskip  Also observe that, similarly to (\ref{eq:N+wm0m}), we have
\begin{equation} \label{eq:N-wm0m}
\N^-_{(w)_m}(x)|_{x=1} = (W^-_\N)_m%({\bf a}; {\bf b}; {\bf c})
 \ \ 
\mathrm{and} \ \  \N^-_{w_0^m}(x)= \widetilde{\N_\tau}(x,y^m).
\end{equation}
We are now ready to prove Theorems 3 and 4. 
It is interesting to note that both results can be proved by the labelled enumeration of weighted $(\G\uparrow \N)^\aut$-structures, with the weights $w$ for Theorem 3 and $w_0$ for Theorem 4.  %We will concentrate on the first result.

\smallskip
\noindent\textbf{Proof of Theorem 3.} 
Consider a $\G\uparrow \N$-structure $(G,G_0)$, with a core $G_0$ in $\G$ and components in $\N$, together with an automorphism $\sigma$. By definition, $\sigma$ should preserve the core $G_0$ and in fact it induces an automorphism $\sigma_0$ of $G_0$. 

We choose to classify these structures according to the cycle type $\underline{n}=(n_1,n_2,\ldots)$ of the core automorphism $\sigma_0$. Also, for enumeration purposes it is preferable to consider the $(\G\uparrow \N)^\aut$-structures with a selected rooting $\rho$ of each cycle of $\sigma_0$. 
Let $(\G\uparrow \N)^\aut_{\underline{n}}$ and $(\G\uparrow \N)^{\aut, \bullet}_{\underline{n}}$ denote the corresponding species.  It is clear that 
$$
(\G\uparrow \N)^{\aut, \bullet}_{\underline{n}}\,(x) 
= 1^{n_1}2^{n_2}\cdots (\G\uparrow \N)^\aut_{\underline{n}}\,(x).
$$
A crucial fact here is that with such a cycle rooting of $\sigma_0$, each edge cycle of $(G_0,\sigma_0)$ inherits a canonical rooting, i.e.~a selection of one of its edges.  
Each of these selected edges can also be canonically oriented, for example according to alphabetical order.  
Then the automorphism $\sigma$ of $G$ will induce a rooted cylindrical $m$-wreath of networks in $\N$ for each rooted cylindrical edge cycle of length $m$ of $\sigma_0$ and a rooted M\"obius $m$-wreath of networks in $\N$ for each rooted M\"obius edge cycle of length $m$. Conversely, the data given by $(G_0,\sigma_0,\rho)$ and the associated family of rooted wreaths completely characterizes the $(\G\uparrow \N)^{\aut, \bullet}_{\underline{n}}$-structure $(G,\sigma,\rho)$. Hence we have
\begin{eqnarray*}
(\G \uparrow \N)_{w}^\aut \, (x) 
& = & \sum_{\underline{n} < \infty} (\G \uparrow \N)_{\underline{n}, {w}}^\aut \,(x)\nonumber \\
& = & \sum_{\underline{n} < \infty} \frac{1}{1^{n_1}2^{n_2}\ldots}(\G \uparrow \N)_{\underline{n}, {w}}^{\aut, \bullet} \,(x) \nonumber \\
& = & \sum_{\underline{n} < \infty} \frac{1}{1^{n_1}2^{n_2}\cdots} \frac{(xa_1)^{n_1}}{n_1!}\frac{(x^2a_2)^{n_2}}{n_2!}\cdots %\nonumber \\& &
\hspace{-2mm} \sum_{G_0\in \mathrm{Fix}\,\G[\underline{n}]} 
\left(\K_1^\bullet(\N)_w^{k_1}\K_2^\bullet(\N)_w^{k_2}\cdots\right) \nonumber\\
& & \ \ \ \ \hspace{58mm} \ \ \ \ \cdot\ %\times\ 
\left(\K_1^{\bullet\mu}(\N)_w^{m_1}\K_2^{\bullet\mu}(\N)_w^{m_2}\ldots %
\right)(x)\nonumber \\
& = & \sum_{\underline{n} < \infty} %\frac{1}{1^{n_1}2^{n_2}\cdots} 
\frac{(xa_1)^{n_1}(x^2a_2)^{n_2}\cdots}{1^{n_1}n_1!2^{n_2}n_2!\cdots}%\frac{}{}
\nonumber \\
& & \hspace{13mm} \sum_{G_0\in \mathrm{Fix}\,\G[\underline{n}]} \N_{(w)_1}^+(x)^{k_1} \N_{{(w)_2}}^+(x^2)^{k_2}\ldots 
\N_{(w)_1}^-(x)^{m_1} \N_{{(w)_2}}^-(x^2)^{m_2}\ldots,
\end{eqnarray*}
where $k_i=k_i(G_0,\underline{n})$ and $m_i=m_i(G_0,\underline{n})$, and, using the representation (\ref{eq:WGter}) of $W_\G$ and the first equations of 
(\ref{eq:N+wm0m}), (\ref{eq:KmNbulletx}), (\ref{eq:KmmuNbulletx}) and (\ref{eq:N-wm0m}), we obtain
\begin{eqnarray}
W_{\G \uparrow \N} ({\bf a}; {\bf b}; {\bf c}) & = & (\G \uparrow \N)_{w}^\aut \ (x)|_{x=1} \nonumber \\
& = & \sum_{\underline{n} < \infty} %\frac{1}{1^{n_1}2^{n_2}\cdots} 
\frac{a_1^{n_1}a_2^{n_2}}{1^{n_1}n_1!2^{n_2}n_2!}\cdots \nonumber \\
& & \hspace{4mm}(\sum_{G_0\in \mathrm{Fix}\,\G[\underline{n}]} \N_{(w)_1}^+(x)^{k_1} \N_{{(w)_2}}^+(x^2)^{k_2}\ldots 
\N_{(w)_1}^-(x)^{m_1} \N_{{(w)_2}}^-(x^2)^{m_2}\ldots)|_{x=1}\nonumber\\
& = & \sum_{\underline{n} < \infty} %\frac{1}{1^{n_1}2^{n_2}\cdots} 
\frac{a_1^{n_1}a_2^{n_2}}{1^{n_1}n_1!2^{n_2}n_2!}\cdots 
\sum_{G_0\in \mathrm{Fix}\,\G[\underline{n}]} (W_\N^+)_1^{m_1}(W_\N^+)_2^{m_2}\ldots 
(W_\N^-)_1^{m_1}(W_\N^-)_2^{m_2}\ldots \nonumber \\
& = & W_\G(a_1,a_2,\ldots; (W_\N^+)_1,(W_\N^+)_2\ldots; (W_\N^-)_1,(W_\N^-)_2\ldots).\nonumber 
\end{eqnarray}
\hfill\rule{2mm}{2mm}

\smallskip
\noindent\textbf{Proof of Theorem 4.}
Although Theorem 4 can be immediately deduced from Theorem 3 by the specialization (\ref{form:8}), a direct proof can also be given, following that of Theorem 3, with the weight $w_0$ instead of $w$. Details are left to the reader. \hfill\rule{2mm}{2mm}
%\smallskip
%\noindent\textbf{Proof of Theorem 4.} Since the representation $\G \uparrow \N$ is canonical, the only possible automorphisms of the graphs in $\G \uparrow \N$ are the automorphisms induced by a core graph $G_0\in \G$. Let $\sigma_0$ be an automorphism of a core graph $G_0$, and the cycle type of $\sigma_0$ is $\underline{n}=(n_1,n_2,\ldots)$. Then we have
%\begin{eqnarray*}
%(\G \uparrow \N)\tilde \ (x,y) & = & (\G \uparrow \N)_{w_0}^\aut \ (x) \nonumber \\
%& = & \sum_{\underline{n} < \infty} (\G \uparrow \N)_{\underline{n}, {w_0}}^\aut \ (x) \nonumber \\
%& = & \sum_{\underline{n} < \infty} \frac{1}{1^{n_1}2^{n_2}\ldots}(\G \uparrow \N)_{\underline{n}, {w_0}}^{\Aut, \bullet} \ (x) \nonumber \\
%& = & \sum_{\underline{n} < \infty} \frac{1}{1^{n_1}2^{n_2}\ldots} \frac{x^{n_1}}{n_1!}\frac{x^{2n_2}}{n_2!}\ldots \nonumber \\
%& & \sum_{G_0\in \mathrm{Fix}\,\G[\underline{n}]} \K_1^\bullet(\N_{w_0})^{m_1}(x)\K_2^\bullet(\N_{w_0})^{m_2}(x)\ldots \nonumber\\
%& & \ \ \ \ \ \ \ \ \ \ \ \ \ \K_1^{\bullet\mu}(\N_{w_0})^{m_1}(x)\K_2^{\bullet\mu}(\N_{w_0})^{m_2}(x)\ldots \nonumber \\
%& = & \sum_{\underline{n} < \infty} \frac{x^{n_1+2n_2+\ldots}}{1^{n_1}n_1!2^{n_2}n_2!\ldots} \nonumber \\
%& & \sum_{G_0\in \mathrm{Fix}\,\G[\underline{n}]} \widetilde{\N}(x,y)^{m_1}\widetilde{\N}(x^2,y^2)^{m_2}\ldots 
%\widetilde{\N}_\tau(x,y)^{m_1}\widetilde{\N}_\tau(x^2,y^2)^{m_2}\ldots \nonumber \\
%& = & W_\G(x,x^2,\ldots; \widetilde{\N}(x,y),\widetilde{\N}(x^2,y^2),\ldots; \widetilde{\N}_\tau(x,y),\widetilde{\N}_\tau(x^2,y^2),\ldots).\nonumber \\
%\end{eqnarray*}
%\hfill\rule{2mm}{2mm}

\bigbreak
{\bf Acknowledgement.} We are thankful to Marie-Eve Harvie for her help in the computations of formulas (\ref{form:NP}) and (\ref{form:NPtau}) and to Timothy Walsh for extending these formulas as well as for numerous helpful discussions.

\end{document}